\documentclass[11pt]{article}
\usepackage{amsfonts,amsmath,amssymb,amscd,amsthm,mathrsfs}
\title{THE SPHERE THEOREMS FOR MANIFOLDS WITH POSITIVE SCALAR CURVATURE\footnote{2010 Mathematics
Subject Classification. 53C20; 53C40.
\newline \indent Keywords: Complete manifolds, sphere
theorem, scalar curvature, Ricci flow, stable current.
\newline\indent Research supported by the NSFC, Grant No. 11071211, 10771187; the Trans-Century Training Programme Foundation for Talents by the
Ministry of Education of China.}}
\author{JUAN-RU GU and HONG-WEI XU}
\date{}
\numberwithin{equation}{section}
\begin{document}
\maketitle
\begin{abstract}
Some new differentiable sphere theorems are obtained via the Ricci
flow and stable currents. We prove that if $M^n$ is a compact
manifold whose normalized scalar curvature and sectional curvature
satisfy the pointwise pinching condition $R_0>\sigma_{n}K_{\max}$,
where $\sigma_n\in (\frac{1}{4},1)$ is an explicit positive
constant, then $M$ is diffeomorphic to a spherical space form. This
gives a partial answer to Yau's conjecture on pinching theorem.
Moreover, we prove that if $M^n(n\geq3)$ is a compact manifold whose
$(n-2)$-th Ricci curvature and normalized scalar curvature satisfy
the pointwise condition $Ric^{(n-2)}_{\min}>\tau_n(n-2)R_0,$ where
$\tau_n\in (\frac{1}{4},1)$ is an explicit positive constant, then
$M$ is diffeomorphic to a spherical space form. We then extend the
sphere theorems above to submanifolds in a Riemannian manifold.
Finally we give a classification of submanifolds with weakly pinched
curvatures, which improves the differentiable pinching theorems due
to Andrews, Baker and the authors.
\end{abstract}

\section{Introduction}
\hspace*{5mm}It plays an important role in global differential
geometry to study curvature and topology of manifolds. The sphere
theorem for Riemannian manifolds was initiated by Rauch \cite{Rauch}
in 1951. During the past six decades, there are many progresses on
sphere theorems for Riemannian manifolds and submanifolds
\cite{Berger,Brendle0,Brendle1,Brendle4,Shiohama,Shiohama1}. The
Brendle-Schoen Differentiable Sphere Theorem
\cite{Brendle2,Brendle3} brought us a big break through in the
investigation of curvature and topology of manifolds. The following
results due to Brendle and Schoen
\cite{Brendle,Brendle3} are very important throughout this paper.\\\\
\textbf{Theorem A(\cite{Brendle}).}\emph{ Let $(M,g_0)$ be a compact
Riemannian manifold of dimension $n(\geq4)$. Assume that
$$R_{1313}+\lambda^2 R_{1414} + R_{2323} + \lambda^2R_{2424}-
2\lambda R_{1234}> 0$$ for all orthonormal four-frames
$\{e_1,e_2,e_3,e_4\}$ and all $\lambda\in[0,1]$. Then the
normalized Ricci flow with initial metric $g_0$
$$\frac{\partial}{\partial t}g(t) = -2Ric_{g(t)} +\frac{2}{n}
r_{g(t)}g(t),$$ exists for all time and converges to a constant
curvature metric as $t\rightarrow\infty$. Here $r_{g(t)}$ denotes
the mean value of the scalar curvature of $g(t)$.}\\\\
\textbf{Theorem B(\cite{Brendle3}).}\emph{ Let $(M,g_0)$ be a
compact, locally irreducible Riemannian manifold of dimension
$n(\geq4)$. Assume that $M\times \mathbb{R}^2$ has nonnegative
isotropic curvature, i.e.,
$$R_{1313}+\lambda^2 R_{1414} +\mu^2 R_{2323} + \lambda^2\mu^2R_{2424}-
2\lambda\mu R_{1234}\geq0$$ for all orthonormal four-frames
$\{e_1,e_2,e_3,e_4\}$ and all $\lambda,\mu\in[-1,1]$. Then
one of the following statements holds:\\
$(i)$ M is diffeomorphic to a spherical space form.\\
$(ii)$ $n=2m$ and the universal cover of M is a K\"{a}hler
manifold biholomorphic
to $\mathbb{C}P^m$.\\
$(iii)$ The universal cover of M is isometric to a compact symmetric space.}\\\\
\hspace*{5mm} On the other hand, some important work on sphere
theorems for manifolds with positive Ricci curvature have been made
by several geometers (see \cite{Berger,Cheeger2,Hamilton,Petersen,
Shen,Shiohama1}, etc.). In 1990's, Cheeger, Colding and Petersen
\cite{Cheeger2,Petersen} proved the following differentiable sphere
theorem for manifolds with positive Ricci curvature.\\\\
\textbf{Theorem C.} \emph{Let $M^n$ be a compact and simply
connected Riemannian $n$-manifold with Ricci curvature $Ric_{M}\ge
n-1.$ Suppose that one of the following conditions
holds:\\
$(i)$ $vol(M)>\omega_n-\varepsilon_1(n)$, where $\omega_n=vol(S^n)$
and $\varepsilon_1(n)$ is some positive constant;\\
$(ii)$ $\lambda_{n+1}<n+\varepsilon_2(n)$, where $\lambda_{n+1}$ is
the $(n+1)$-th eigenvalue of $M$ and $\varepsilon_2(n)$
is some positive constant.\\
Then $M$ is diffeomorphic to $S^n$.}\\\\
\hspace*{5mm}Let $K(\pi)$ be the sectional curvature of $M$ for
2-plane $\pi\subset T_xM$, $Ric(u)$ the Ricci curvature of $M$ for
unit vector $u\in U_xM$. Set $K_{\max}(x):=\max_{\pi\subset
T_{x}M}K(\pi)$, $Ric_{\min}(x):=\min_{u\in U_{x}M}Ric(u).$ Inspired
by Shen's topological sphere theorem \cite{Shen}, the authors
\cite{Xu3} obtained the following differentiable sphere theorem for
manifolds of positive Ricci curvatures.\\\\
\textbf{Theorem D.} \emph{Let $M^n$ be a compact Riemannian
$n$-manifold. If $Ric_{\min}>\delta_{n}(n-1)K_{\max}$, where
$\delta_{n}=1-\frac{6}{5(n-1)}$,
 then $M$ is diffeomorphic to a spherical space form. In particular, if M is simply connected, then M is diffeomorphic
to $S^n$ .}\\\\
\hspace*{5mm}Let $M^{n}$ be a submanifold in a Riemannian manifold
$\overline{M}^{N}$. Denote by $H$ and $S$ the mean curvature and the
squared length of the second fundamental form  of $M$, respectively.
Denote by $\overline{K}(\pi)$ the sectional curvature of
$\overline{M}$ for 2-plane $\pi(\subset T_x\overline{M})$. Set
$\overline{K}_{\max}(x):=\max_{\pi\subset
T_{x}\overline{M}}\overline{K}(\pi)$,
$\overline{K}_{\min}(x):=\min_{\pi\subset
T_{x}\overline{M}}\overline{K}(\pi)$. In \cite{Xu4}, Xu and Zhao
obtained some differentiable sphere theorems for complete
submanifolds in higher codimensions via the Ricci flow and stable
currents. Recently the authors \cite{Xu2} proved the
following differentiable sphere theorem for complete submanifolds with strictly pinched curvatures.\\\\
\textbf{Theorem E.} \emph{Let $M^n$ be an $n$-dimensional complete
submanifold in an $N$-dimensional Riemannian manifold
$\overline{M}^{N}$.   If
$S<\frac{8}{3}\Big(\overline{K}_{\min}-\frac{1}{4}\overline{K}_{\max}\Big)+
\frac{n^2H^{2}}{n-1}$,
 then $M$ is diffeomorphic to
a spherical space form or $\mathbb{R}^n$. In particular, if M is
simply connected, then M is diffeomorphic
to $S^n$ or $\mathbb{R}^n$.}\\\\
\hspace*{5mm}The purpose of this paper is to prove some new
differentiable sphere theorems for Riemannian manifolds and
submanifolds. In
Section 3, we prove the following differentiable sphere theorem  for compact manifolds with positive scalar curvature. \\\\
\textbf{Theorem 1.1.} \emph{Let $M^n$ be an $n(\geq 3)$-dimensional
compact Riemannian manifold. Denote by $R_0$  the
normalized scalar curvature of $M$. Assume that one of the following pointwise conditions holds: \\
$(i)$ $R_0>\sigma_n K_{\max};$\\
 $(ii)$ $K_{\min}>\eta_n R_0.$\\
 Then $M$ is diffeomorphic to a spherical space form. In particular,
if M is simply connected, then M is diffeomorphic to $S^n$. Here
\begin{eqnarray*}&&\sigma_n=\left\{\begin{array}{ll}1-\frac{4}{n(n-1)} \mbox{\ \ \ \ for\ \ } n=3,
    \\
1-\frac{12}{5n(n-1)} \mbox{\ \ \ for\ \ } n\geq 4,
\end{array}\right. \\
 &&\eta_n=1-\frac{6}{n^2-n+6}.\end{eqnarray*}
} \\
\hspace*{5mm}Theorem 1.1 improves Theorem D and gives a partial
answer to Yau's conjecture on pointwise pinching theorem (See
\cite{Yau0}, Problem
12). Moreover, we obtain the following theorem.\\\\
\textbf{Theorem 1.2.} \emph{Let $M^n$ be an $n(\geq 3)$-dimensional
compact Riemannian manifold. Denote by $R_0$ and $Ric^{(n-2)}$ the
normalized scalar curvature and the
$(n-2)$-th Ricci curvature of $M$. Assume that one of the following pointwise conditions holds: \\
 $(i)$ $(n-2)R_0>\mu_{n}Ric_{\max}^{(n-2)}$;\\
 $(ii)$ $Ric^{(n-2)}_{\min}>\tau_n(n-2)R_0 .$\\
Then $M$ is diffeomorphic to a spherical space form. In particular,
if M is simply connected, then M is diffeomorphic to $S^n$. Here
\begin{eqnarray*}
&&\mu_n=1-\frac{6}{n(n-1)(n+1)},\\
&&\tau_n=\max\{1-\frac{12}{(n-2)(5n^2-11n-6)},\,0\}.\end{eqnarray*}
} \\
\textbf{Remark 1.1.} Note that differentiable structures of Einstein
manifolds are very rich. If pinching conditions (i) and (ii) in
Theorem 1.2 are replaced by $(n-1)R_0>\tilde{\mu}_n Ric_{\max}$ and
$Ric_{\min}>\tilde{\tau}_n(n-1)R_0 $ respectively, where
$\tilde{\mu}_n, \tilde{\tau}_n\in (\frac{1}{4},1)$, it's impossible
to obtain the same assertion. Therefore, the pinching conditions in
Theorem 1.2 are the weakest in this sense.\\\\
 \hspace*{5mm}In Section 4 we extend the sphere theorems above to
submanifolds in a Riemannian manifold with arbitrary codimension
(Theorems 4.4 and 4.5). In Section 5, we obtain a differentiable
sphere
theorem for submanifolds with weakly pinched curvatures, stated as: \\\\
\textbf{Theorem 1.3.} \emph{Let $M^n$ be an $n(\geq3)$-dimensional
compact submanifold in an $N$-dimensional Riemannian manifold
$\overline{M}^{N}$. Assume that $M$ satisfies one of the following conditions:\\
$(i)$ $\overline{K}_{\min}(x_{0})-\frac{1}{4}\overline{K}_{\max}(x_{0})\neq 0$ for some point $x_{0}\in M $, and $S\leq\frac{8}{3}\Big(\overline{K}_{\min}-\frac{1}{4}\overline{K}_{\max}\Big)+
\frac{n^2H^{2}}{n-1}$;\\
$(ii)$ $\overline{K}_{\min}(x)-\frac{1}{4}\overline{K}_{\max}(x)=0$ for any point $x\in M $, $S\leq \frac{n^2H^{2}}{n-1}$ and
the strict inequality holds for some point $x_{0}\in M$.\\
Then $M$ is diffeomorphic
a spherical space form. In particular, if M is simply connected, then M is diffeomorphic to $S^n$.}\\\\
\hspace*{5mm}Furthermore, we prove the following classification theorem of submanifolds with weakly pinched curvatures in space forms.\\\\
\textbf{Theorem 1.4.} \emph{Let $M^n$ be an $n(\geq4)$-dimensional
oriented complete submanifold in an $N$-dimensional simply connected
space form $F^{N}(c)$ with $c\ge0$. Assume that its scalar curvature
$R\ge
(n+1)(n-2)c+\frac{n^2(n-2)}{n-1}H^{2}$, where $c+H^2>0$. We have\\
$(i)$ If $c=0$, then M is either diffeomorphic $S^n$, $\mathbb{R}^n$, or locally isometric to $S^{n-1}(r)\times \mathbb{R}$.\\
$(ii)$ If M is compact, then M is diffeomorphic to $S^n$ .}\\\\
\textbf{Remark 1.2.} The pinching condition in Theorem 1.4 is
equivalent to $S\leq 2c+\frac{n^2H^2}{n-1}$. Theorems 1.3 and 1.4
improve the differentiable pinching theorems due to Andrews-Baker
and the authors \cite{Andrews,Xu2}. \\\\
\hspace*{5mm}It should be mentioned that the second author
introduced the results above in his invited talk at the Fifth
International Congress of Chinese Mathematicians held in Beijing from December 17 to December 22, 2010.\\\\

\section{Notation and lemmas}
\hspace*{5mm}Let $M^{n}$ be an $n$-dimensional submanifold in an
$N$-dimensional Riemannian manifold $\overline{M}^{N}$. We shall make use
of the following convention on the range of indices.
$$ 1\leq A,B,C,\ldots\leq N;\ 1\leq i,j,k,\ldots\leq n;$$  $$\mbox{if\, } N\geq n+1,\, \ n+1\leq
\alpha,\beta,\gamma,\ldots\leq N .$$ For an arbitrary fixed point
$x\in M\subset \overline{M}$, we choose an orthonormal local frame field
$\{e_{A}\}$ in $\overline{M}^{N}$ such that $e_{i}$'s are tangent to $M$.
Denote by $\{\omega_{A}\}$ the dual frame field of $\{e_{A}\}$. Let
$$Rm=\sum_{i,j,k,l}{R_{ijkl}}\omega_{i}\otimes\omega_{j}\otimes\omega_{k}\otimes\omega_{l},\hspace*{13mm}$$
$$\overline{Rm}=\sum_{A,B,C,D}\overline{R}_{ABCD}\omega_{A}\otimes\omega_{B}\otimes\omega_{C}\otimes\omega_{D},$$
be the Riemannian curvature tensors of $M$ and $\overline{M}$,
respectively. Denote by $h$ and $\xi$ the second fundamental form
and the mean curvature vector of $M$. When $N=n$, $h$ and $\xi$ are
identically equal to zero. When $N\ge n+1$, we set
$$h=\sum_{\alpha,i,j}h^{\alpha}_{ij}\omega_{i}\otimes\omega_{j}\otimes
e_{\alpha},  \ \ \ \
 \xi=\frac{1}{n}\sum_{\alpha,i}h^{\alpha}_{ii}e_{\alpha}.$$The
squared norm $S$ of the second fundamental form and the mean
curvature $H$ of $M$ are given by $
S:=\sum_{\alpha,i,j}(h^{\alpha}_{ij})^{2},\,\,H:=|\xi|.$ Then we
have the Gauss equation
 \begin{equation}R_{ijkl}=\overline{R}_{ijkl}+<h(e_i,e_k),h(e_j,e_l)>-<h(e_i,e_l),h(e_j,e_k)>.
\end{equation}Denote by $K(\cdot)$, $\overline{K}(\cdot)$,
$Ric(\cdot)$, $\overline{Ric}(\cdot)$, $R$ and $\overline{R}$ the
sectional curvatures, the Ricci curvatures and the scalar curvatures
of $M$ and $\overline{M}$, respectively. Then we have
$$Ric(e_{i})=\sum_{j}R_{ijij},\, \, \,  \overline{Ric}(e_{A})=\sum_{B}
\overline{R}_{ABAB},$$
$$R=\sum_{i,j}R_{ijij}, \, \, \,   \overline{R}=\sum_{A,B}
\overline{R}_{ABAB}.$$
Set $K_{\min}(x)=\min_{\pi\subset
T_{x}M}K(\pi)$, $K_{\max}(x)=\max_{\pi\subset T_{x}M}K(\pi)$,
  $\overline{K}_{\min}(x)=\min_{\pi\subset
T_{x}\overline{M}}\overline{K}(\pi)$,
$\overline{K}_{\max}(x)=\max_{\pi\subset
T_{x}\overline{M}}\overline{K}(\pi)$.
  Then by Berger's inequality(See e.g.
\cite{Brendle0}, Proposition 1.9), we
have\begin{equation}|R_{ijkl}|\leq\frac{2}{3}(K_{\max}-K_{\min})\end{equation}
 for all distinct indices $i$, $j$, $k$, $l$, and  \begin{equation}|\overline{R}_{ABCD}|\leq\frac{2}{3}(\overline{K}_{\max}-\overline{K}_{\min})\end{equation}
 for all distinct indices $A$, $B$, $C$, $D$. We set
\begin{eqnarray*}
&&Ric_{\min}(x)=\min_{u\in U_{x}M}Ric(u),\,\,  \overline{Ric}_{\min}(x)=\min_{u\in U_{x}\overline{M}} \overline{Ric}(u),\\
&&Ric_{\max}(x)=\max_{u\in U_{x}M}Ric(u),\,\,
\overline{Ric}_{\max}(x)=\max_{u\in U_{x}\overline{M}}
\overline{Ric}(u).\end{eqnarray*} For any unit tangent vector $u\in
U_{x}M$ at point $x\in M,$ let $V_{x}^k$ be a $k$-dimensional
subspace of $T_{x}M$ satisfying $u\perp V_{x}^k$. Choose an
orthonormal basis $\{e_{i}\}$ in $T_{x}M$ such that $e_{j_0}=u,\,\,
span\{e_{j_1},\ldots,e_{j_k}\}=V_{x}^k$, where the indices $1\le
j_0,j_1,\ldots,j_k\le n$ are distinct with each other. We set
 \begin{eqnarray}&&Ric^{(k)}(u;V_{x}^k)=Ric ^{(k)}([e_{j_0},\ldots,e_{j_{k}}])=\sum_ {q=1}^{k}R_{j_0j_{q}j_0j_{q}},\nonumber\\&&Ric^{(k)}_{\min}(x)=\min_{u\in U_{x}M}\min_{u\perp
V_{x}^k\subset T_{x}M}
Ric^{(k)}(u;V_{x}^k),\nonumber\\&&Ric^{(k)}_{\max}(x)=\max_{u\in
U_{x}M}\max_{u\perp V_{x}^k\subset T_{x}M}Ric^{(k)}(u;V_{x}^k).\end{eqnarray}
 We extend  an orthonormal
$s$-frame $\{e_{j_0},\ldots,e_{j_{s-1}}\}$ in $T_{x}M$ to
$(k+1)$-frame $\{e_{j_0},\ldots,e_{j_{k}}\}$ for $1\leq s\leq
k+1\leq n$ and set
\begin{eqnarray}&&R^{(k,s)}([e_{j_0},\ldots,e_{j_{k}}])=\sum_ {p=0}^{s-1}\sum_ {q=0}^{k}R_{j_pj_{q}j_pj_{q}},\nonumber\\
&&R^{(k,s)}_{\min}(x)=\min_{
 \{e_{j_0},\ldots,e_{j_{k}}\}\subset T_{x}M} R^{(k,s)}([e_{j_0},\ldots,e_{j_{k}}]),\nonumber\\
 &&R^{(k,s)}_{\max}(x)=\max_{
 \{e_{j_0},\ldots,e_{j_{k}}\}\subset T_{x}M} R^{(k,s)}([e_{j_0},\ldots,e_{j_{k}}]).\end{eqnarray}
 \begin{eqnarray}&&Ric^{[s]}([e_{j_0},\ldots,e_{j_{n-1}}])=R^{(n-1,s)}([e_{j_0},\ldots,e_{j_{n-1}}])=\sum_ {p=0}^{s-1}\sum_ {q=0}^{n-1}R_{j_pj_{q}j_pj_{q}},\nonumber\\
&&Ric^{[s]}_{\min}(x)=\min_{
 \{e_{j_0},\ldots,e_{j_{n-1}}\}\subset T_{x}M} Ric^{[s]}([e_{j_0},\ldots,e_{j_{n-1}}]),\nonumber\\
 &&Ric^{[s]}_{\max}(x)=\max_{
 \{e_{j_0},\ldots,e_{j_{n-1}}\}\subset T_{x}M} Ric^{[s]}([e_{j_0},\ldots,e_{j_{n-1}}]).\end{eqnarray}
\begin{eqnarray}&&R^{(k)}([e_{j_0},\ldots,e_{j_{k}}])=R^{(k,k+1)}([e_{j_0},\ldots,e_{j_{k}}])=\sum_ {p=0}^{k}\sum_ {q=0}^{k}R_{j_pj_{q}j_pj_{q}},\nonumber\\
 &&R^{(k)}_{\min}(x)=\min_{
 \{e_{j_0},\ldots,e_{j_{k}}\}\subset T_{x}M} R^{(k)}([e_{j_0},\ldots,e_{j_{k}}]) ,\nonumber\\
&&R^{(k)}_{\max}(x)=\max_{
 \{e_{j_0},\ldots,e_{j_{k}}\}\subset T_{x}M}  R^{(k)}([e_{j_0},\ldots,e_{j_{k}}]) .\end{eqnarray}
 \textbf{Definition 2.1.} \emph{We call $Ric^{(k)}(u;V_{x}^k)$, $R^{(k,s)}([e_{j_0},\ldots,e_{j_{k}}])$, $Ric^{[s]}([e_{j_0},\ldots,e_{j_{n-1}}])$
 and $R^{(k)}([e_{j_0},\ldots,e_{j_{k}}])$ the $k$-th Ricci curvature, $(k,s)$-curvature, $s$-th weak Ricci curvature and $k$-th scalar curvatrure of $M$, respectively.}\\\\
The geometry and topology of $k$-th Ricci curvature was initiated by
Hartman \cite{Hartman} in 1979, and developed by Wu \cite{Wu} and
Shen \cite{Shen,Shen1}, etc.. By the definition above, it is seen
that the Ricci curvature of $M$ is equal to the $(n-1)$-th Ricci
curvature, $(n-1,1)$-curvature and $1$-th weak Ricci curvature; the
scalar curvature of $M$ is equal to $(n-1,n)$-curvature, $n$-th weak
Ricci curvature and $(n-1)$-th scalar curvature. For any unit
tangent vector $u\in U_{x}\overline{M}$ at point $x\in
\overline{M},$ let $V_{x}^k$ be a $k$-dimensional subspace of
$T_{x}\overline{M}$ satisfying $u\perp V_{x}^k$. Choose an
orthonormal basis $\{e_{A}\}$ in $T_{x}\overline{M}$ such that
$e_{A_0}=u,\,\, span\{e_{A_1},\ldots,e_{A_k}\}=V_{x}^k$, where the
indices $1\le A_0,A_1,\ldots,A_k\le N$ are distinct with each other.
We define the $k$-th Ricci curvature  as follows.
\begin{equation}\overline{Ric}^{(k)}(u;V_{x}^k)=\sum_ {q=1}^{k}\overline{R}_{A_0A_{q}A_0A_{q}}.\end{equation} We extend an
orthonormal $s$-frame $\{e_{A_0},\ldots,e_{A_{s-1}}\}$ in
$T_{x}\overline{M}$ to $(k+1)$-frame $\{e_{A_0},\ldots,e_{A_{k}}\}$
for $1\leq s\leq k+1\leq N$ and defined the $(k,s)$-curvature,
$s$-th weak Ricci curvature and $k$-th scalar curvature of
$\overline{M}$ as follows.
\begin{eqnarray}&&\overline{R}^{(k,s)}([e_{A_0},\ldots,e_{A_{k}}])=\sum_ {p=0}^{s-1}\sum_ {q=0}^{k}\overline{R}_{A_pA_{q}A_pA_{q}}\nonumber\\
&&\overline{Ric}^{[s]}([e_{A_0},\ldots,e_{A_{N-1}}])=\overline{R}^{(N-1,s)}([e_{A_0},\ldots,e_{A_{N-1}}])=\sum_ {p=0}^{s-1}\sum_ {q=0}^{N-1}\overline{R}_{A_pA_{q}A_pA_{q}}\nonumber\\
&&\overline{R}^{(k)}([e_{A_0},\ldots,e_{A_{k}}])=\overline{R}^{(k,k+1)}([e_{A_0},\ldots,e_{A_{k}}])=\sum_ {p=0}^{k}\sum_ {q=0}^{k}\overline{R}_{A_pA_{q}A_pA_{q}}.\end{eqnarray}
Denote by $\overline{Ric}^{(k)}_{\min}(x), \overline{R}^{(k,s)}_{\min}(x),\overline{Ric}^{[s]}_{\min}(x),
\overline{R}^{(k)}_{\min}(x)$ and $\overline{Ric}^{(k)}_{\max}(x), \overline{R}^{(k,s)}_{\max}(x),
 \overline{Ric}^{[s]}_{\max}(x),
\overline{R}^{(k)}_{\max}(x)$ the minimum and maximum of the curvatures defined above at point $x\in \overline{M}$.\\\\
\hspace*{5mm}We choose an orthonormal frame
$\{e_{1},e_{2},\cdots,e_{n}\}$ such that $u=e_{n}$ and
$Ric^{(k)}(u;V^k_x)=\sum_{i=1}^{k}R_{inin}$, where
$V^k_x=span\{e_1,e_2,\cdots,e_k\},\,\,
 1\leq k\leq n-1$. In particular, we see that
$Ric^{(n-1)}(u;V^{n-1}_x)=Ric(u)$ and
$Ric^{(1)}(u;V^1_x)=K(\pi)$, where $\pi=span\{e_1,e_n\}$. Then we have the following lemma.\\\\
\textbf{Lemma 2.1.} \emph{Let $M^{n}$ be an $n$-dimensional complete
 submanifold in an
$N$-dimensional Euclidean space $\mathbb{R}^{N}$. If $S\leq\frac{n^2H^{2}}{n-1}$, $H\neq 0$, then\\
$(i)(\cite{Shiohama2, Xu2})$ $Ric^{(k)}(u;V^k_x)\geq0.$\\
$(ii)$ For each point x $\in$ M there exists a unit vector u such
that $Ric^{(k)}(u;V^k_x)=0$  for some integer $k\in[2,n-1]$ if and
only if $H$ is a constant and M is isometric
to $S^{n-1}\Big(\frac{n-1}{nH}\Big)\times \mathbb{R} $.}\\\\
\textbf{Proof.}  If $k=1$, the assertion follows from the result
in \cite{Xu2}. Now we discuss the case for $2\leq k\leq n-1$. Choose
an orthonormal frame $\{e_{1},e_{2},\cdots,e_{N}\}$ such that
$e_{n+1}$ is parallel to the
 mean curvature vector $\xi$. Then
 \begin{eqnarray}
n^2H^2&=&\Big(\sum_{i=1}^nh_{ii}^{n+1}\Big)^2\nonumber\\
&=&(n-1)\Big[\sum_{i=1}^{n}(h^{n+1}_{ii})^2+\sum_{i\neq
j}(h_{ij}^{n+1})^2+\sum_{\alpha=n+2}^{N}\sum_{i,
j=1}^{n}(h_{ij}^{\alpha})^2+\frac{n^2H^2}{n-1}-S\Big].\end{eqnarray}
Note that for $l\neq n$
\begin{eqnarray*}
\Big(\sum_{i=1}^nh_{ii}^{n+1}\Big)^2&\leq&(n-1)\Big[(h^{n+1}_{ll}+h^{n+1}_{nn})^2+\sum_{i\neq l,n}(h^{n+1}_{ii})^2\Big]\\
&=&(n-1)\Big[\sum_{i=1}^{n}(h^{n+1}_{ii})^2+2h^{n+1}_{ll}h^{n+1}_{nn}\Big].\end{eqnarray*}
This together with (2.10) implies
\begin{equation}
2h^{n+1}_{ll}h^{n+1}_{nn}\geq\sum_{i\neq
j}(h_{ij}^{n+1})^2+\sum_{\alpha=n+2}^{N}\sum_{i,
j=1}^{n}(h_{ij}^{\alpha})^2+\frac{n^2H^2}{n-1}-S \end{equation} for
$l\neq n$. The equality holds if and only if
$h^{n+1}_{ii}=h^{n+1}_{ll}+h^{n+1}_{nn}$ for $i\neq l,n$. This
together with (2.1) implies\begin{eqnarray*}
Ric^{(k)}(u;V^k_x)&=&\sum_{i=1}^{k}R_{inin}=\sum_{i=1}^k\sum_{\alpha=n+1}^{N}[h^{\alpha}_{ii}h^{\alpha}_{nn}-(h^{\alpha}_{in})^2]\\
 &\geq&k\sum_{\alpha=n+1}^{N}\sum_{1\leq i<j<n}(h_{ij}^{\alpha})^2+(k-1)\sum_{\alpha=n+1}^{N}\sum_{i=1}^{n-1}
(h_{in}^{\alpha})^2\\
 &&+\frac{k-1}{2}\sum_{\alpha=n+2}^{N}\sum_{i=1}^{n-1}
(h_{ii}^{\alpha})^2+\frac{1}{2}\sum_{\alpha=n+2}^{N}\sum_{i=1}^{k}(h^{\alpha}_{ii}+h^{\alpha}_{nn})^2+\frac{k}{2}\Big(\frac{n^2H^2}{n-1}-S\Big)\\
 &\geq&\frac{k}{2}\Big(\frac{n^2H^2}{n-1}-S\Big) . \end{eqnarray*}
 If $S\leq\frac{n^2H^{2}}{n-1}$, then $Ric^{(k)}(u;V^k_x)\geq 0$. The equality holds if and only if
 $$h^{\alpha}_{ij}=0,\,1\le i,j\le n,\, \alpha\neq n+1;\,\,\, h^{n+1}_{ij}=0,\, i\neq j,\,1\le i,j\le n;$$   $$h^{n+1}_{nn}=0;\,\,\, h^{n+1}_{ii}=\frac{nH}{n-1},\, 1\le i\le n-1.$$
 Hence $M$ has essential codimension one. Since the shape operator of $M$ has one eigenvalue of multiplicity $n-1$ and the other eigenvalue is zero,
 it follows from a result due to Deprez(See \cite{Deprez}, Corollary) that $H$ is a constant and $M$ is isometric to $S^{n-1}\Big(\frac{n-1}{nH}\Big)\times \mathbb{R} .$
This completes the proof.\\\\
\hspace*{5mm}The following nonexistence theorem for stable currents
in a compact Riemannian manifold $M$ isometrically immersed into
$F^{N}(c)$ is employed to eliminate the homology groups
$H_{q}(M;\mathbb{Z})$ for $0<q<n$, which was initiated by
Lawson-Simons \cite{Lawson2} and extended by Xin \cite{Xin}.\\\\
\textbf{Theorem 2.1.} \emph{Let $M^{n}$ be a compact submanifold in
$F^{N}(c)$ with $c\geq 0$. Assume that$$
\sum_{k=q+1}^{n}\sum_{i=1}^{q}[2|h(e_{i},e_{k})|^{2}-\langle
h(e_{i},e_{i}),h(e_{k},e_{k})\rangle]<q(n-q)c$$ holds for any
orthonormal basis $\{e_{i}\}$ of $T_xM$ at any point $x\in M$,
where q is an integer satisfying $0<q<n$. Then there does not exist
any stable q-currents. Moreover,
$$H_{q}(M;\mathbb{Z})=H_{n-q}(M;\mathbb{Z})=0,$$ where $H_{i}(M;\mathbb{Z})$ is the $i$-th homology group of M with integer
coefficients, and $\pi_{1}(M)=0$ when $q=1$.}\\\\
 \hspace*{5mm}From the proof of Lemma 2 in \cite{Shiohama2}, we have
 \begin{eqnarray}
&&\sum_{k=q+1}^{n}\sum_{i=1}^{q}[2|h(e_{i},e_{k})|^{2}-\langle
h(e_{i},e_{i}),h(e_{k},e_{k})\rangle] \nonumber\\
&\leq&\frac{q(n-q)}{n}\bigg[S-2nH^{2}+\frac{\sqrt{n}|2q-n|}{\sqrt{q(n-q)}}H\sqrt{S-nH^{2}}\bigg]\nonumber\\
&\leq&\frac{q(n-q)}{n}\bigg[S-2nH^{2}+\frac{n(n-4)}{\sqrt{2n(n-2)}}H\sqrt{S-nH^{2}}\bigg],
\end{eqnarray}
for $n\geq 4$ and $1<q<n-1$.
This together with Theorem 2.1 implies the following.\\\\
\textbf{Lemma 2.2.} \emph{Let $M^{n}$ be an $n(\geq4)$-dimensional
 compact submanifold in an Euclidean space $\mathbb{R}^{N}$.
 If $S<\frac{n^2H^2}{n-2}$, then $$H_{q}(M;\mathbb{Z})=0,\mbox{\ \ for all \ \
 }1<q<n-1.\vspace*{-3mm}$$}\\
 \textbf{Lemma 2.3(\cite{Harish}).} \emph{Let $M$ be a compact Riemannian
manifold of dimension n. If
M has nonnegative isotropic curvature and has positive isotropic curvature for some point in M,
then M admits a metric with positive isotropic curvature.}\\\\
 \textbf{Lemma 2.4.} \emph{Let $M$ be a compact Riemannian
manifold of dimension n. If
$M\times \mathbb{R}^2$ has nonnegative isotropic curvature, and if $M$ has positive Ricci curvature and
isotropic curvature, then M is diffeomorphic to a spherical space form.}\\\\
\textbf{Proof.} By the assume that $M$ has positive Ricci curvature,
the universal cover $\widetilde{M}$ of $M$ is compact. Since $M$ has positive isotropic
curvature, $\widetilde{M}$ also has positive isotropic curvature.
 Note that $\widetilde{M}$ is simply connected. It follows from a theorem due to Micallef and
 Moore \cite{Micallef} that $\widetilde{M}$ is homeomorphic to $S^n$.
 Therefore, $M$ is locally irreducible and the symmetric metric of $\widetilde{M}$ would have to be of
 positive constant curvature. Moreover, when $n$ is even, a theorem due to Micallef and Wang \cite{Micallef2} states that if $\widetilde{M}$ has positive isotropic curvature, then
 $H^2(\widetilde{M},\mathbb{R})=0$.
Hence $\widetilde{M}$ can not be a K\"{a}hler manifold. This together with Theorem B implies that $M$ is diffeomorphic to a spherical space form.\\\\
 \section{Manifolds of positive scalar curvature}
 \hspace*{5mm} In this section, we will give the proof of Theorem
 1.1. More generally, we will prove Theorem 3.3.
We first prove the following lemma for compact manifolds.\\\\
\textbf{Lemma 3.1.} \emph{Let $M^n$ be an $n(\geq4)$-dimensional
compact Riemannian manifold. Denote by $R^{(k)}(\cdot)$ and
$R^{(k,s)}(\cdot)$ the $k$-th scalar curvature and $(k,s)$-curvature of $M$. If one of the following conditions holds:\\
$(i)$ $R^{(k)}_{\min}>\Big(k^2+k-\frac{24}{7}\Big)K_{\max}$ for
some integer $k\in[3, n-1],$\\
$(ii)$ $R^{(k,s)}_{\min}>\frac{s(7k^2+7k-24)}{7(k+1)}K_{\max}$ for some integers $k\in[3, n-1]$ and $s\in[2,k+1],$\\
then $\pi_{k}(M)=0$ for $2 \leq k\leq [\frac{n}{2}]$. In particular, if M is
simply connected, then M is homeomorphic to a sphere.}\\\\
 \textbf{Proof.}
(i) It follows from (2.7) that $$R^{(k)}_{\min}\leq
2K_{\min}+[k(k+1)-2]K_{\max}.$$ Then we have
\begin{equation}K_{\min}\geq\frac{1}{2}[ R^{(k)}_{\min}-(k^2+k-2)K_{\max}].\end{equation}
 Suppose $\{e_1,e_2,e_3,e_4\}$ is an
orthonormal four-frame.
From (2.2), (3.1) and the assumption we get
\begin{eqnarray}&&R_{1313}+ R_{1414}+R_{2323}+ R_{2424}-2R_{1234}\nonumber\\
 &\geq& \frac{1}{2}\{ R^{(k)}_{\min}-[k(k+1)-8] K_{\max}\}-\frac{4}{3}(K_{\max}-K_{\min})\nonumber\\
 &\geq& \frac{1}{2}\{R^{(k)}_{\min}-[k(k+1)-8]K_{\max}\}-\frac{2}{3}[k(k+1) K_{\max}-R_{\min}^{(k)}]\nonumber\\
&\geq&\frac{7}{6}\Big[ R^{(k)}_{\min}-\Big(k^2+k-\frac{24}{7}\Big) K_{\max}\Big] \nonumber\\
&>&0.\end{eqnarray}Hence $M$ has positive isotropic curvature. By a
result due to Micallef and Moore \cite{Micallef}, we have
$\pi_{k}(M)=0$ for $2\leq k\leq[\frac{n}{2}]$. In particular, if $M$
is simply connected, then $M$ is homeomorphic
to a sphere.\\
\hspace*{5mm}(ii) By Definition 2.1, we have
\begin{equation}\frac{R^{(k)}_{\min}}{k(k+1)}\geq
\frac{R^{(k,s)}_{\min}}{ks}.
\end{equation}
This together with the assumption implies
\begin{equation}R^{(k)}_{\min}\geq
\frac{k+1}{s}R^{(k,s)}_{\min}>\Big(k^2+k-\frac{24}{7}\Big)K_{\max}.\end{equation}
Then the assertion follows from (i).\\
\hspace*{5mm}This completes the proof of Lemma 3.1.\\\\
\hspace*{5mm} By taking $k=n-1$ in Lemma 3.1, we have the following.\\\\
\textbf{Theorem 3.1.} \emph{Let $M^n$ be an $n(\geq4)$-dimensional
compact Riemannian manifold. Denote by $Ric^{[s]}(\cdot)$ and $R_0(\cdot)$ the $s$-th weak Ricci curvature and normalized scalar curvature of M. If one of the following conditions holds:\\
 $(i)$ $Ric^{[s]}_{\min}>\frac{s
(7n^2-7n-24)}{7n}K_{\max}$ for some integer  $s\in[2, n],$\\
$(ii)$ $R_0>\Big[1-\frac{24}{7n(n-1)}\Big]K_{\max},$\\
 then $\pi_{k}(M)=0$ for $2 \leq k\leq [\frac{n}{2}]$. In particular, if M is simply connected, then M is homeomorphic to a sphere.}\\\\
 \textbf{Corollary 3.1.} \emph{Let $M^n$ be an $n$-dimensional
compact and simply connected Riemannian manifold, where $4\leq n\leq
6$. Denote by $R_0$ the normalized scalar curvature of M. If
$R_0>\Big[1-\frac{24}{7n(n-1)}\Big]K_{\max}$,
 then M is diffeomorphic to $S^n$.}\\\\
\textbf{ Proof.} It follows from Theorem 3.1 that $M$ has positive isotropic curvature. A theorem due to Hamilton \cite{Hamilton2} says that a 4-dimensional compact simply connected manifold with positive isotropic curvature is diffeomorphic to $S^4$.
It is well known that there is only one differentiable structure on $S^n$, $n=5,6$. This together with Theorem 3.1 implies $M$ is diffeomorphic to $S^n$ for $n=5,6$. This proves the corollary.\\\\
 \textbf{Lemma 3.2.} \emph{Let $M^n$ be an $n(\geq 4)$-dimensional
compact Riemannian manifold.  Denote by $R^{(k)}(\cdot)$ and
$R^{(k,s)}(\cdot)$  the $k$-th scalar curvature and $(k,s)$-curvature
 of $M$, respectively. If one of the following conditions
holds:\\ $(i)$ $R^{(k)}_{\min}>\Big(k^2+k-\frac{12}{5}\Big)K_{\max}$ for some integer $k\in[2, n-1]$,\\
$(ii)$ $R^{(k,s)}_{\min}>\frac{s
(5k^2+5k-12)}{5(k+1)}K_{\max}$ for some integers $k\in[2, n-1]$
and $s\in[2,k+1],$\\
 then $M$ is diffeomorphic to a spherical space form.}\\\\
\textbf{Proof.} (i)
 Suppose $\{e_1,e_2,e_3,e_4\}$ is an
orthonormal four-frame and $\lambda\in \mathbb{R}$.
From (2.2) and (3.1) we obtain
\begin{eqnarray}&&R_{1313}+ R_{2323}-|R_{1234}|\nonumber\\
&\geq&
 \frac{1}{2}\Big[R^{(k)}_{\min}-2\sum_{i<j\neq 3 }^{k+1} R_{ijij}]-\frac{2}{3}( K_{\max}- K _{\min}) \nonumber\\
 &\geq&\frac{1}{2}[R^{(k)}_{\min}-(k^2+k-4)K_{\max}]-\frac{1}{3}[k(k+1)K_{\max}- R^{(k)}_{\min}]\nonumber\\
 &\geq&\frac{5}{6}\Big[ R^{(k)}_{\min}-\Big(k^2+k-\frac{12}{5}\Big) K_{\max}\Big].\end{eqnarray}
 Using the same argument as above, we get
\begin{equation}R_{1414}+ R_{2424}-|R_{1234}|\geq\frac{5}{6}\Big[ R^{(k)}_{\min}-\Big(k^2+k-\frac{12}{5}\Big) K_{\max}\Big] .\end{equation}
 From (3.5), (3.6) and the assumption, we have
\begin{eqnarray}&&R_{1313}+\lambda^2 R_{1414} + R_{2323} +
\lambda^2R_{2424}- 2\lambda R_{1234}\nonumber\\
&\geq&R_{1313}+ R_{2323}-|R_{1234}|+\lambda^2(R_{1414}+ R_{2424}-|R_{1234}|)\nonumber\\
&\geq&\frac{5(1+\lambda^2)}{6}\Big[ R^{(k)}_{\min}-\Big(k^2+k-\frac{12}{5}\Big) K_{\max}\Big]\nonumber\\
&>&0.
\end{eqnarray}
This together with Theorem A implies that $M$ is diffeomorphic to a spherical space form.\\
\hspace*{5mm}(ii) From (3.3) and the assumption we know that \begin{equation}R^{(k)}_{\min}\geq \frac{k+1}{s}R^{(k,s)}_{\min}>\Big(k^2+k-\frac{12}{5}\Big)K_{\max}.\end{equation}
Hence the conclusion follows from (i).\\
\hspace*{5mm}This proves Lemma 3.2.\\\\
\textbf{Lemma 3.3.} \emph{Let $M^n$ be an $n(\geq4)$-dimensional
compact Riemannian manifold. If its $k$-th Ricci curvature satisfies one of the following conditions:\\
 $(i)$ $Ric^{(k)}_{\min}>\frac{5k-6}{5k-1}Ric^{(k+1)}_{\max};$ \\
$(ii)$ $Ric^{(k)}_{\min}>\frac{(5k-6)(k+1)}{(5k-1)k}Ric^{(k)}_{\max}$,\\
 where k is some integer in $[2,n-2]$, then $M$ is diffeomorphic to a spherical space form.}\\\\
\textbf{Proof.} (i) From (2.4), we obtain \begin{equation}K_{\max}\leq Ric^{(k+1)}_{\max}-Ric^{(k)}_{\min},\end{equation}
and $$Ric^{(k)}_{\min}\leq
K_{\min}+(k-1)K_{\max},$$ which implies
\begin{eqnarray}K_{\min}&\geq &Ric^{(k)}_{\min}-(k-1)K_{\max}\nonumber\\
&\geq & kRic^{(k)}_{\min}-(k-1)Ric^{(k+1)}_{\max}.\end{eqnarray}
Suppose $\{e_1,e_2,e_3,e_4\}$ is an
orthonormal four-frame and $\lambda\in \mathbb{R}$.
Then we have from (2.2), (2.4), (3.9) and (3.10) that
\begin{eqnarray}&&R_{1313}+ R_{2323}-|R_{1234}|\nonumber\\
&\geq&
Ric^{(k)}_{\min}-\sum_{i= 3}^{k+1}R_{i3i3}-\frac{2}{3}(K_{\max}-K_{\min})\nonumber\\
&\geq&Ric^{(k)}_{\min}-\Big(k-\frac{4}{3}\Big)[Ric^{(k+1)}_{\max}-Ric^{(k)}_{\min}]+\frac{2}{3}[kRic^{(k)}_{\min}-(k-1)Ric^{(k+1)}_{\max}]\nonumber\\
&\geq&\frac{5k-1}{3}\Big[Ric^{(k)}_{\min}-\frac{5k-6}{5k-1}Ric^{(k+1)}_{\max}\Big].\end{eqnarray}
 Similarly, we get
\begin{equation}R_{1414}+ R_{2424}-|R_{1234}|\geq\frac{5k-1}{3}\Big[Ric^{(k)}_{\min}-\frac{5k-6}{5k-1}Ric^{(k+1)}_{\max}\Big].\end{equation}
From the (3.11), (3.12) and the assumption we obtain
\begin{eqnarray}&&R_{1313}+\lambda^2 R_{1414} + R_{2323} +
\lambda^2R_{2424}- 2\lambda R_{1234}\nonumber\\
&\geq&R_{1313}+ R_{2323}-|R_{1234}|+\lambda^2(R_{1414}+ R_{2424}-|R_{1234}|)\nonumber\\
&\geq&\frac{(1+\lambda^2)(5k-1)}{3}\Big[ Ric^{(k)}_{\min}-\frac{5k-6}{5k-1}Ric^{(k+1)}_{\max}\Big]\nonumber\\
&>&0.\end{eqnarray}
This together with Theorem A implies that $M$ is diffeomorphic a spherical space form.\\
\hspace*{5mm}(ii) From (2.4) we have \begin{equation}
\frac{Ric^{(k)}_{\max}}{k}\geq \frac{Ric^{(k+1)}_{\max}}{k+1},\end{equation}
which together with the assumption implies $$ 0<Ric^{(k)}_{\min}-\frac{(5k-6)(k+1)}{(5k-1)k}Ric^{(k)}_{\max}\leq Ric^{(k)}_{\min}-\frac{5k-6}{5k-1}Ric^{(k+1)}_{\max}.$$
Hence the assertion follows from (i).\\
\hspace*{5mm}This proves Lemma 3.3.\\\\
\hspace*{5mm}Taking $k=n-2$ in condition (i) of Lemma 3.3, we have the following theorem.\\\\
\textbf{Theorem 3.2.} \emph{Let $M^n$ be an $n(\geq4)$-dimensional
compact Riemannian manifold. If its $(n-2)$-th Ricci curvature and Ricci
curvature satisfy
$$Ric^{(n-2)}_{\min}>\frac{5n-16}{5n-11}Ric_{\max},$$
 then $M$ is diffeomorphic to a spherical space form. In particular, if M is simply connected, then M is diffeomorphic
to $S^n$ .}\\\\
\textbf{Lemma 3.4.} \emph{Let $M^n$ be an $n(\geq4)$-dimensional
compact Riemannian manifold. Denote by $Ric^{(k)}(\cdot),$
$R^{(k)}(\cdot)$ and
$R^{(k,s)}(\cdot)$ the $k$-th Ricci curvature, $k$-th scalar curvature
and $(k,s)$-curvature of $M$, respectively. If one of the following conditions holds:\\
$(i)$ $Ric^{(k )}_{\min}>\frac{5k-6}{5k^2+9k-8}R^{(k+1)}_{\max}$
for some integer $k\in[2,n-2]$,\\
 $(ii)$ $Ric^{(k)}_{\min}>\frac{(k+2)(5k-6)}{s(5k^2+9k-8)}R^{(k+1,s)}_{\max}$
for some integers $k\in[2,n-2]$ and $s\in[2,k+2]$,\\
then $M$ is diffeomorphic to a spherical space form.}\\\\
\textbf{Proof.} (i) It follows from (2.4) and (2.7) that
\begin{equation}K_{\max}\leq
\frac{1}{2}[R^{(k+1)}_{\max}-(k+3)Ric^{(k)}_{\min}],\end{equation}
and $$Ric^{(k)}_{\min}\leq K_{\min}+(k-1)K_{\max}.$$  Then we have
\begin{eqnarray}K_{\min}&\geq &Ric^{(k)}_{\min}-(k-1)K_{\max}\nonumber\\
&\geq& \frac{1}{2}[(k^2+2k-1)Ric^{(k)}_{\min}-(k-1)R^{(k+1)}_{\max}]
.\end{eqnarray} Suppose $\{e_1,e_2,e_3,e_4\}$ is an orthonormal
four-frame and $\lambda\in \mathbb{R}$. Combing (2.2), (2.4), (3.15)
and (3.16), we have
\begin{eqnarray}&&R_{1313}+ R_{2323}-|R_{1234}|\nonumber\\
&\geq&
Ric^{(k)}_{\min}-\sum_{i= 3}^{k+1}R_{i3i3}-\frac{2}{3}(K_{\max}-K_{\min})\nonumber\\
&\geq&Ric^{(k)}_{\min}-\Big(\frac{k}{2}-\frac{2}{3}\Big)[R^{(k+1)}_{\max}-(k+3)Ric^{(k)}_{\min}]
\nonumber\\&&+\frac{1}{3}[(k^2+2k-1)Ric^{(k)}_{\min}-(k-1)R^{(k+1)}_{\max}]\nonumber\\
&\geq&\frac{5k^2+9k-8}{6}\Big[Ric^{(k)}_{\min}-\frac{5k-6}{5k^2+9k-8}R^{(k+1)}_{\max}\Big].\end{eqnarray}
 By a similar argument, we obtain
\begin{equation}R_{1414}+ R_{2424}-|R_{1234}|\geq\frac{5k^2+9k-8}{6}\Big[Ric^{(k)}_{\min}-\frac{5k-6}{5k^2+9k-8}R^{(k)}_{\max}\Big].\end{equation}
From (3.17), (3.18) and the assumption we obtain
\begin{eqnarray}&&R_{1313}+\lambda^2 R_{1414} + R_{2323} +
\lambda^2R_{2424}- 2\lambda R_{1234}\nonumber\\
&\geq&R_{1313}+ R_{2323}-|R_{1234}|+\lambda^2(R_{1414}+ R_{2424}-|R_{1234}|)\nonumber\\
&\geq&\frac{(1+\lambda^2)(5k^2+9k-8)}{6}\Big[Ric^{(k)}_{\min}-\frac{5k-6}{5k^2+9k-8}R^{(k)}_{\max}\Big]\nonumber\\
&>&0.\end{eqnarray}
This together with Theorem A implies that $M$ is diffeomorphic a spherical space form.\\
 \hspace*{5mm}(ii) We get from (2.5) and (2.7) that \begin{equation} \frac{R^{(k+1,s)}_{\max}}{s(k+1)}\geq\frac{R^{(k+1 )}_{\max}}{(k+1)(k+2)},\end{equation}
which together with the assumption implies
\begin{eqnarray}Ric^{(k)}_{\min}&>&\frac{(k+2)(5k-6)}{s(5k^2+9k-8)}R^{(k+1,s)}_{\max}\nonumber\\
&\geq&\frac{5k-6}{5k^2+9k-8}R^{(k+1 )}_{\max}. \end{eqnarray}
The assertion follows from (i).\\
\hspace*{5mm}This proves the lemma.\\\\
 \textbf{Theorem 3.3.} \emph{Let $M^n$ be an $n(\geq4)$-dimensional
compact Riemannian manifold. Denote by $Ric^{[s]}(\cdot),$
$Ric^{(k)}(\cdot)$ and $K(\cdot)$ the $s$-th weak Ricci curvature, $k$-th Ricci curvature
and sectional curvature of $M$, respectively. Suppose one of the following conditions holds:\\
$(i)$ $Ric^{[s]}_{\min}>\frac{s(5n^2-5n-12)}{5n}K_{\max}$ for some integer $s\in[2,n]$;\\
$(ii)$ $Ric^{[s]}_{\min}>\frac{s(n^2+2n+3)}{n(n+1)}Ric_{\max}^{(n-2)}$ for some integer $s\in[2,n]$;\\
$(iii)$ $K_{\min}>\frac{1}{s(n-1)+6}Ric_{\max}^{[s]}$ for some integer $s\in[2,n]$;\\
$(iv)$ $Ric^{(n-2)}_{\min}>\frac{n(5n-16)}{s(5n^2-11n-6)}Ric^{[s]}_{\max}$ for some integer $s\in[2,n]$.\\
Then the normalized Ricci flow with initial metric $g_0$
$$\frac{\partial}{\partial t}g(t) = -2Ric_{g(t)} +\frac{2}{n}
r_{g(t)}g(t),$$ exists for all time and converges to a constant
curvature metric as $t\rightarrow\infty$. Moreover, $M$ is diffeomorphic to a spherical space form. In particular, if M is simply connected, then M is diffeomorphic
to $S^n$.}\\\\
\textbf{ Proof.}
(i) Taking $k=n-1$ in condition (ii) of Lemma 3.2, we get the conclusion.\\
\hspace*{5mm}(ii) Since \begin{eqnarray} K_{\min}&\geq &\frac{1}{2}[R-(n+1)Ric^{(n-2)}_{\max}]\nonumber\\
&\geq&\frac{1}{2}\Big[\frac{nRic_{\min}^{[s]}}{s}-(n+1)Ric^{(n-2)}_{\max}\Big],\end{eqnarray}
we obtain\begin{eqnarray}K_{\max}&\leq& Ric^{(n-2)}_{\max}-(n-3)K_{\min}\nonumber\\
&\leq&\frac{n^2-2n-1}{2}Ric_{\max}^{(n-2)}-\frac{n(n-3)}{2s}Ric_{\min}^{[s]}\end{eqnarray}
Suppose $\{e_1,e_2,e_3,e_4\}$ is an orthonormal four-frame and
$\lambda\in \mathbb{R}$. It follows from (2.2), (3.22) and (3.23)
that
\begin{eqnarray}&&R_{1313}+ R_{2323}-|R_{1234}|\nonumber\\
&\geq&
2K_{\min}-\frac{2}{3}(K_{\max}-K_{\min})\nonumber\\
&\geq&\frac{4}{3}\Big[\frac{nRic_{\min}^{[s]}}{s}-(n+1)Ric^{(n-2)}_{\max}\Big]-\frac{2}{3}\Big[\frac{n^2-2n-1}{2}Ric_{\max}^{(n-2)}-\frac{n(n-3)}{2s}Ric_{\min}^{[s]}\Big]\nonumber\\
&\geq&\frac{1}{3}\Big[\frac{n(n+1)}{s}Ric_{\min}^{[s]}-(n^2+2n+3)Ric_{\max}^{(n-2)}\Big].\end{eqnarray}
 Similarly, we get
\begin{equation}R_{1414}+ R_{2424}-|R_{1234}|\geq\frac{1}{3}\Big[\frac{n(n+1)}{s}Ric_{\min}^{[s]}-(n^2+2n+3)Ric_{\max}^{(n-2)}\Big].\end{equation}
 From (3.24), (3.25) and the assumption we obtain
\begin{eqnarray}&&R_{1313}+\lambda^2 R_{1414} + R_{2323} +
\lambda^2R_{2424}- 2\lambda R_{1234}\nonumber\\
&\geq&R_{1313}+ R_{2323}-|R_{1234}|+\lambda^2(R_{1414}+ R_{2424}-|R_{1234}|)\nonumber\\
&\geq&\frac{1+\lambda^2}{3}\Big[\frac{n(n+1)}{s}Ric_{\min}^{[s]}-(n^2+2n+3)Ric_{\max}^{(n-2)}\Big]\nonumber\\
&>&0.\end{eqnarray}
This together with Theorem A implies $M$ is diffeomorphic to a spherical space form. \\
\hspace*{5mm}(iii) By Definition 2.1, we get
\begin{equation}K_{\max}\leq
\frac{1}{2}\{Ric_{\max}^{[s]}-[s(n-1)-2]K_{\min}\}.\end{equation}Suppose
$\{e_1,e_2,e_3,e_4\}$ is an orthonormal four-frame and $\lambda\in
\mathbb{R}$. It follows from (2.2) and (3.27) that
\begin{eqnarray}&&R_{1313}+ R_{2323}-|R_{1234}|\nonumber\\
&\geq&
2K_{\min}-\frac{2}{3}(K_{\max}-K_{\min})\nonumber\\
&\geq&\frac{8}{3}K_{\min}-\frac{1}{3}[Ric_{\max}^{[s]}-(sn-s-2)K_{\min}]\nonumber\\
&\geq& \frac{1}{3}[(sn-s+6)K_{\min}-Ric_{\max}^{[s]}].\end{eqnarray}
 A similar discussion implies that
\begin{equation}R_{1414}+ R_{2424}-|R_{1234}|\geq\frac{1}{3}[(sn-s+6)K_{\min}-Ric_{\max}^{[s]}].\end{equation}
 From (3.28), (3.29) and the assumption we obtain
\begin{eqnarray}&&R_{1313}+\lambda^2 R_{1414} + R_{2323} +
\lambda^2R_{2424}- 2\lambda R_{1234}\nonumber\\
&\geq&R_{1313}+ R_{2323}-|R_{1234}|+\lambda^2(R_{1414}+ R_{2424}-|R_{1234}|)\nonumber\\
&\geq&\frac{1+\lambda^2}{3}[(sn-s+6)K_{\min}-Ric_{\max}^{[s]}]\nonumber\\
&>&0.\end{eqnarray}
This together with Theorem A implies that $M$ is diffeomorphic a spherical space form.\\
\hspace*{5mm}(iv) The assertion follows by taking $k=n-2$ in (ii) of Lemma 3.4.\\
\hspace*{5mm}This proves the theorem.\\\\
\hspace*{5mm}We are now in a position to prove Theorem 1.1.\\\\
\textbf{Proof of
Theorem 1.1.} (i) If $n=3$, for any unit tangent vector $u\in
U_xM$ at $x\in M$, we choose an orthonormal three-frame
$\{e_1,e_2,e_3\}$ such that $e_3=u$. Then from the assumption we
obtain
\begin{eqnarray*}Ric(u)&=&R_{1313}+ R_{2323} \\
&=&
 \frac{1}{2}(R-2R_{1212}) \\
&\geq& \frac{1}{2}(R-2K_{\max}) \\
&>& 0.\end{eqnarray*} This together with Hamilton's theorem
\cite{Hamilton} implies that $M$ is diffeomorphic to a spherical
space form. When $n\geq4$,  the assertion follows by taking $k=n-1$
in (i) of Lemma 3.2.\\
\hspace*{5mm}(ii) If $n=3$, the assertion follows from Hamilton's work \cite{Hamilton}. Thus from now on
we assume that $n\geq 4$. By taking $s=n$ in (iii) of Theorem 3.3, we conclude that $M$ is diffeomorphic a spherical space form.\\
\hspace*{5mm}This proves Theorem 1.1.\\\\
\hspace*{5mm}In 1990, Yau \cite{Yau0} proposed the following conjecture (see also \cite{Schoen2,Yau1}).\\\\
\textbf{Yau Conjecture I.} \emph{Let $M^n$ be a compact and simply
connected Riemannian manifold. Denote by $R_0$ the normalized scalar
curvature of M. If $K_{\min}>\frac{n-1}{n+2}R_0$,
then $M$ is diffeomorphic to $S^n$. }\\\\
\hspace*{5mm} If $n=2,3$, the answer is affirmative. If the pinching
constant in Yau Conjecture I is replaced by
$\eta_n=\frac{n^2-n}{n^2-n+6}$, Theorem 1.1 gives an affirmative
answer. The following example shows that $\frac{n-1}{n+2}$ is the
best possible pinching constant for the conjecture in even
dimensions($\ge4$).\\\\
\textbf{Example 3.1.} Let $R_0$ be the normalized scalar curvature
of a Riemannian manifold. By a direct computation, we have the
normalized scalar curvatures of the compact rank one symmetric
spaces (CROSS) with standard metrics.
\begin{eqnarray*}&& R_0(\mathbb{C}P^m)=\frac{m+1}{4m-2},\, \,\dim_\mathbb{R}(\mathbb{C}P^m)=2m,\,m\ge2; \\
&&
R_0(\mathbb{H}P^m)=\frac{m+2}{4m-1},\,\,\dim_\mathbb{R}(\mathbb{H}P^m)=4m,\,m\ge2;
\\
&&R_0(\mathbb{O}P^2)=\frac{3}{5},\,\,\dim_\mathbb{R}(\mathbb{O}P^2)=16.
\end{eqnarray*}
On the other hand,
$$K_{\min}(\mathbb{C}P^m)=K_{\min}(\mathbb{H}P^m)=K_{\min}(\mathbb{O}P^2)=\frac{1}{4},$$
and these are not homeomorphic to $S^n$. Therefore,
$\frac{n-1}{n+2}$ is the best possible pinching constant for Yau Conjecture I in even dimensions($\ge4$).\\\\
\textbf{Yau Conjecture II.} \emph{Let $M^n$ be a compact and simply
connected Riemannian manifold. Denote by $R_0$ the normalized scalar
curvature of M. If $K_{M}\ge \frac{n-1}{n+2}$ and $R_0\le 1$, then
$M$ is
either diffeomorphic to $S^n$, or isometric to the complex projective space $\mathbb{C}\emph{P}^m$ with $n=2m$. }\\\\
\hspace*{5mm} Recently the authors \cite{Xu1} proved the following
optimal
rigidity theorem for Einstein manifolds, which provides an evidence for Yau Conjectures I and II.\\\\
 \textbf{Theorem 3.4.} \emph{Let $M$ be an $n(\geq 4)$-dimensional compact
Einstein manifold with normalized scalar curvature $R_0:=c$. If
$K_{\min}\ge \frac{n-1}{n+2}R_0>0$, then $M$ is locally symmetric.
In particular, if $M$ is simply connected, then $M$ is isometric to
either the standard n-sphere $S^n(\frac{1}{\sqrt{c}})$
or the complex projective space $\mathbb{C}\emph{P}^m(c)$.}\\\\
\hspace*{5mm} Motivated by Theorems 1.1 and Example 3.1, we would
like to propose the
following conjectures.\\\\
\textbf{Conjecture A.} \emph{Let $M^n(n\geq4)$ be a compact
Riemannian manifold. If $R_0>\frac{3}{5}K_{\max}$,
then $M$ is diffeomorphic to a spherical space form. In particular, if M is simply connected, then M is diffeomorphic to $S^n$.}\\\\
\textbf{Conjecture B.} \emph{Let $M^n(n\geq4)$ be an even
dimensional compact and
simply connected Riemannian manifold.
If $K_M\leq 1$ and $R_0\geq c_n$, where
$$c_n=\left\{\begin{array}{ll}
\frac{n+2}{4(n-1)}\mbox{\ \ \ for $n=4$  or  $4k+2,$ \, $k\in
\mathbb{Z^+}$,}
    \\
\frac{n+8}{4(n-1)}\mbox{\ \ \ for $n=4k$, \, $k\in
\mathbb{Z^+}\bigcap[2,\infty)$ and $k\ne4$,}
    \\
\frac{3}{5}\mbox{\, \, \, \, \, \, for $n=16,$}
\end{array}\right. $$
then $M$ is either diffeomorphic to $S^n$, or isometric to a compact
rank one symmetric space.}\\\\
\textbf{Proof of
Theorem 1.2.} (i) If $n=3$, the assertion follows from (i) of Theorem 1.1. If $n\geq 4$, the conclusion follows from (ii) of Theorem 3.3 by taking $s=n$.\\
 \hspace*{5mm}(ii) If $n=3$, it follows from Hamilton's work \cite{Hamilton}. If $n\geq 4$, by taking $k=n-2$ in (i) of Lemma 3.4, we get the conclusion.\\
 \hspace*{5mm}This completes the proof of Theorem 1.2.\\\\
\section{Sphere theorems for compact submanifolds}
\hspace*{5mm} In this section, we extend the sphere theorems in
Section 3 to submanifolds in Riemannian manifolds with arbitrary
codimension. For compact submanifolds, we prove the following lemma.\\\\
\textbf{Lemma 4.1.} \emph{Let $M^n$ be an $n(\geq4)$-dimensional
compact submanifold in an $N$-dimensional Riemannian manifold
$\overline{M}^{N}$. Assume that $M$ satisfies one of the following conditions:\\
$(i)\, S<\frac{7}{6}\Big[\overline{R}^{(k)}_{\min}-\Big(k^2+k-\frac{24}{7}\Big)\overline{K}_{\max}\Big]+
\frac{n^2H^{2}}{n-2}\mbox{\, for some integer\, }k\in[3, N-1];$\\
$(ii)\, S<\frac{7(k+1)}{6s}\Big[\overline{R}^{(k,s)}_{\min}-\frac{s(7k^2+7k-24)}{7(k+1)}\overline{K}_{\max}\Big]+
\frac{n^2H^{2}}{n-2}$ for some integers $k\in[3, N-1]$ and  $s\in[2, k+1].$\\
 Then $\pi_{k}(M)=0$ for  $2\leq k\leq[\frac{n}{2}]$. In particular, if M is simply connected, then M is homeomorphic
to a sphere.}\\\\
\textbf{Proof.} (i) Since $$\overline{R}^{(k)}_{\min}\leq
2\overline{K}_{\min}+[k(k+1)-2]\overline{K}_{\max},$$ we have
\begin{equation}\overline{K}_{\min}\geq\frac{1}{2}[\overline{R}^{(k)}_{\min}-(k^2+k-2)\overline{K}_{\max}].\end{equation} Setting $S_{\alpha}= \sum_{i,j=1}^{n}(h^{\alpha}_{ij})^{2},$ we
know that
\begin{equation}
 \Big(\sum_{i=1}^nh_{ii}^{\alpha}\Big)^2=(n-2)\Big[\sum_{i=1}^{n}(h^{\alpha}_{ii})^2+\sum_{i\neq
j}(h_{ij}^{\alpha})^2+\frac{(\sum_{i=1}^nh_{ii}^{\alpha})^2}{n-2}-S_{\alpha}\Big].\end{equation}
Note that for all distinct $p, q , m, l$
\begin{eqnarray*}
\Big(\sum_{i=1}^nh_{ii}^{\alpha}\Big)^2&\leq&(n-2)\Big[(h^{\alpha}_{pp}+h^{\alpha}_{qq})^2+(h^{\alpha}_{mm}+h^{\alpha}_{ll})^2+\sum_{i\neq p,q,m,l}(h^{\alpha}_{ii})^2\Big]\\
&=&(n-2)\Big[\sum_{i=1}^{n}(h^{\alpha}_{ii})^2+2h^{\alpha}_{pp}h^{\alpha}_{qq}+2h^{\alpha}_{mm}h^{\alpha}_{ll}\Big].\end{eqnarray*}
This together with (4.2) implies
\begin{equation}
2h^{\alpha}_{pp}h^{\alpha}_{qq}+2h^{\alpha}_{mm}h^{\alpha}_{ll}\geq\sum_{i\neq
j}(h_{ij}^{\alpha})^2+\frac{(\sum_{i=1}^nh_{ii}^{\alpha})^2}{n-2}-S_{\alpha},\end{equation}
for all distinct $p,q,m,l$.
 Suppose $\{e_1,e_2,e_3,e_4\}$ is an
orthonormal four-frame and $\lambda\in \mathbb{R}$.
From (2.1), (2.3), (4.1), (4.3) and the assumption we get
\begin{eqnarray}&&R_{1313}+ R_{1414}+R_{2323}+ R_{2424}-2R_{1234}\nonumber\\
&=&\overline{R}_{1313}+ \overline{R}_{1414}+\overline{R}_{2323}+ \overline{R}_{2424}-2\overline{R}_{1234}
 +\sum_{\alpha}\Big[h_{11}^{\alpha}h_{33}^{\alpha}+h_{22}^{\alpha}h_{44}^{\alpha}
 +h_{22}^{\alpha}h_{33}^{\alpha}\nonumber\\
 &&+h_{11}^{\alpha}h_{44}^{\alpha}
 - (h_{13}^{\alpha})^2
-(h_{23}^{\alpha})^2-(h^{\alpha}_{24})^2-(h^{\alpha}_{14})^2-2(h^{\alpha}_{13}h^{\alpha}_{24}-h^{\alpha}_{14}h^{\alpha}_{23})\Big]\nonumber\\
&\geq& \frac{1}{2}\{\overline{R}^{(k)}_{\min}-[k(k+1)-8]\overline{K}_{\max}\}-\frac{4}{3}(\overline{K}_{\max}-\overline{K}_{\min})\nonumber\\
&&+\sum_{\alpha}\Big[\sum_{i\neq
j}(h_{ij}^{\alpha})^2+\frac{(\sum_{i=1}^nh_{ii}^{\alpha})^2}{n-2}-S_{\alpha}
-2(h_{13}^{\alpha})^2
-2(h_{23}^{\alpha})^2-2(h^{\alpha}_{24})^2-2(h^{\alpha}_{14})^2\Big]\nonumber\\
&\geq& \frac{1}{2}[\overline{R}^{(k)}_{\min}-(k^2+k-8)\overline{K}_{\max}] -\frac{2}{3}[k(k+1)\overline{K}_{\max}-\overline{R}^{(k)}_{\min}]+ \frac{n^2H^2}{n-2}-S\nonumber\\
&\geq&\frac{7}{6}\Big[\overline{R}^{(k)}_{\min}-\Big(k^2+k-\frac{24}{7}\Big)\overline{K}_{\max}\Big]+ \frac{n^2H^2}{n-2}-S\nonumber\\
&>&0.\end{eqnarray}
 Therefore $M$ has positive isotropic curvature. From Micallef and Moore's theorem
 \cite{Micallef}, we get
$\pi_{k}(M)=0$ for $2\leq k\leq[\frac{n}{2}]$. In particular, if $M$
is simply connected, then $M$ is homeomorphic
to a sphere. \\
\hspace*{5mm}(ii) Notice that
\begin{equation}\frac{\overline{R}^{(k)}_{\min}}{k(k+1)}\geq
\frac{\overline{R}^{(k,s)}_{\min}}{ks}.
\end{equation}
We have
\begin{equation}\overline{R}^{(k)}_{\min}-\Big(k^2+k-\frac{24}{7}\Big)\overline{K}_{\max}\geq
\frac{k+1}{s}\overline{R}^{(k,s)}_{\min}-\Big(k^2+k-\frac{24}{7}\Big)\overline{K}_{\max}.\end{equation}
The assertion follows from (i), (4.6) and the assumption.\\
\hspace*{5mm}This proves Lemma 4.1.\\\\
\hspace*{5mm}By taking $k=N-1$ in Lemma 4.1, we get the following theorem.\\\\
 \textbf{Theorem 4.1.} \emph{Let $M^n$ be an $n(\geq4)$-dimensional
compact submanifold in an $N$-dimensional Riemannian manifold
$\overline{M}^{N}$. Assume that $M$ satisfies one of the following conditions:\\
$(i)\, S<\frac{7}{6}\Big[\overline{R} -\Big(N^2-N-\frac{24}{7}\Big)\overline{K}_{\max}\Big]+
\frac{n^2H^{2}}{n-2}$;\\
$(ii)\, S<\frac{7N}{6s}\Big[\overline{Ric}^{[s]}_{\min}-\frac{s(7N^2-7N-24)}{7N}\overline{K}_{\max}\Big]+
\frac{n^2H^{2}}{n-2}$ for some integer $s\in[2, N].$\\
 Then $\pi_{k}(M)=0$ for  $2\leq k\leq[\frac{n}{2}]$. In particular, if M is simply connected, then M is homeomorphic
to a sphere.}\\\\
\textbf{Lemma 4.2.} \emph{Let $M^n$ be an $n(\geq4)$-dimensional
compact submanifold in an $N$-dimensional Riemannian manifold
$\overline{M}^{N}$. Suppose that $M$ satisfies one of the following conditions:\\
$(i)$ $S<\frac{5}{6}\Big[\overline{R}^{(k)}_{\min}-\Big(k^2+k-\frac{12}{5}\Big)\overline{K}_{\max}\Big]+
\frac{n^2H^{2}}{n-1},$ for some integer $k\in[2,N-1]$;\\
$(ii)$ $S<\frac{5(k+1)}{6s}\Big[\overline{R}^{(k,s)}_{\min}-\frac{s( 5k^2+5k-12)}{5(k+1)}\overline{K}_{\max}\Big] +
\frac{n^2H^{2}}{n-1},$
 for some integers $k\in[2,N-1]$ and $s\in[2,k+1]$.\\
 Then $M$ is diffeomorphic to a spherical space form. }\\\\
\textbf{Proof.}
(i) Setting $S_{\alpha}=\sum_{i,j=1}^{n}(h^{\alpha}_{ij})^{2},$ we
have
\begin{equation}
 \Big(\sum_{i=1}^nh_{ii}^{\alpha}\Big)^2=(n-1)\Big[\sum_{i=1}^{n}(h^{\alpha}_{ii})^2+\sum_{i\neq
j}(h_{ij}^{\alpha})^2+\frac{(\sum_{i=1}^nh_{ii}^{\alpha})^2}{n-1}-S_{\alpha}\Big].\end{equation}
Note that for $m\neq l$
\begin{eqnarray*}
\Big(\sum_{i=1}^nh_{ii}^{\alpha}\Big)^2&\leq&(n-1)\Big[(h^{\alpha}_{mm}+h^{\alpha}_{ll})^2+\sum_{i\neq m,l}(h^{\alpha}_{ii})^2\Big]\\
&=&(n-1)\Big[\sum_{i=1}^{n}(h^{\alpha}_{ii})^2+2h^{\alpha}_{mm}h^{\alpha}_{ll}\Big].\end{eqnarray*}
This together with (4.7) implies
\begin{equation}
2h^{\alpha}_{mm}h^{\alpha}_{ll}\geq\sum_{i\neq
j}(h_{ij}^{\alpha})^2+\frac{(\sum_{i=1}^nh_{ii}^{\alpha})^2}{n-1}-S_{\alpha}\end{equation}
 for all distinct $m, l$, and the equality holds if and only if $h^{\alpha}_{ii}=h^{\alpha}_{mm}+h^{\alpha}_{ll}$ for all $i\neq m,l$.
 Suppose $\{e_1,e_2,e_3,e_4\}$ is an
orthonormal four-frame and $\lambda\in \mathbb{R}$.
From (2.1), (2.3), (4.1) and (4.8), we obtain
\begin{eqnarray}&&R_{1313}+ R_{2323}-|R_{1234}|\nonumber\\
&\geq&
 \frac{1}{2}\Big[\overline{R}^{(k)}_{\min}-2\sum_{A<B\neq 3 }^{k+1}\overline{R}_{ABAB}]-\frac{2}{3}(\overline{K}_{\max}-\overline{K}_{\min})+\sum_{\alpha}\Big[h_{11}^{\alpha}h_{33}^{\alpha}+h_{22}^{\alpha}h_{33}^{\alpha}\nonumber\\
&&-\frac{3}{2}(h_{13}^{\alpha})^2
-\frac{3}{2}(h_{23}^{\alpha})^2-\frac{1}{2}(h^{\alpha}_{24})^2-\frac{1}{2}(h^{\alpha}_{14})^2\Big]\nonumber\\
&\geq&\frac{1}{2}\Big[\overline{R}^{(k)}_{\min}-(k^2+k-4)\overline{K}_{\max}]-\frac{1}{3}[k(k+1)\overline{K}_{\max}-\overline{R}^{(k)}_{\min}]\nonumber\\
&&+\sum_{\alpha}\Big[\sum_{i\neq
j}(h_{ij}^{\alpha})^2+\frac{(\sum_{i=1}^nh_{ii}^{\alpha})^2}{n-1}-S_{\alpha}
-\frac{3}{2}(h_{13}^{\alpha})^2
-\frac{3}{2}(h_{23}^{\alpha})^2-\frac{1}{2}(h^{\alpha}_{24})^2-\frac{1}{2}(h^{\alpha}_{14})^2\Big]\nonumber\\
&\geq&\frac{5}{6}\Big[\overline{R}^{(k)}_{\min}-\Big(k^2+k-\frac{12}{5}\Big)\overline{K}_{\max}\Big]+ \frac{n^2H^2}{n-1}-S.\end{eqnarray}
 Similarly, we get
\begin{equation}R_{1414}+ R_{2424}-|R_{1234}|\geq\frac{5}{6}\Big[\overline{R}^{(k)}_{\min}-\Big(k^2+k-\frac{12}{5}\Big)\overline{K}_{\max}\Big]+
\frac{n^2H^2}{n-1}-S.\end{equation}
 From (4.9) and (4.10), we obtain
\begin{eqnarray}&&R_{1313}+\lambda^2 R_{1414} + R_{2323} +
\lambda^2R_{2424}- 2\lambda R_{1234}\nonumber\\
&\geq&R_{1313}+ R_{2323}-|R_{1234}|+\lambda^2(R_{1414}+ R_{2424}-|R_{1234}|)\nonumber\\
&\geq&(1+\lambda^2)\Big\{\frac{5}{6}\Big[\overline{R}^{(k)}_{\min}-\Big(k^2+k-\frac{12}{5}\Big)\overline{K}_{\max}\Big]+
\frac{n^2H^2}{n-1}-S\Big\} .
\end{eqnarray}
This together with Theorem A and the assumption implies $M$ is diffeomorphic to a spherical space form.  \\
 \hspace*{5mm}(ii)
  From (4.5) we have  \begin{equation}\overline{R}^{(k)}_{\min}-\Big(k^2+k-\frac{12}{5}\Big)\overline{K}_{\max}\geq \frac{k+1}{s}\overline{R}^{(k,s)}_{\min}-\Big(k^2+k-\frac{12}{5}\Big)\overline{K}_{\max}.\end{equation}
Therefore the assertion follows from (i), (4.12) and the assumption.\\
\hspace*{5mm}This completes the proof.\\\\
 \textbf{Lemma 4.3.} \emph{Let $M^n$ be an $n(\geq4)$-dimensional
compact submanifold in an $N$-dimensional Riemannian manifold
$\overline{M}^{N}$. Suppose that $M$ satisfies one of the following conditions:\\
$(i)$
$S<\frac{5k-1}{3}\Big[\overline{Ric}^{(k)}_{\min}-\frac{5k-6}{5k-1}\overline{Ric}^{(k+1)}_{\max}\Big]+
\frac{n^2H^{2}}{n-1}$ for some integer $k\in[2, N-2]$;\\
$(ii)$ $S<\frac{5k-1}{3}\Big[\overline{Ric}^{(k)}_{\min}-\frac{(5k-6)(k+1)}{(5k-1)k}\overline{Ric}^{(k)}_{\max}\Big]+
\frac{n^2H^{2}}{n-1}$, for some integer $k\in[2, N-2]$.\\ Then $M$
is diffeomorphic to a spherical space form. }\\\\
\textbf{Proof.}
Since \begin{equation}\overline{K}_{\max}\leq \overline{Ric}^{(k+1)}_{\max}-\overline{Ric}^{(k)}_{\min},\end{equation}
and $$\overline{Ric}^{(k)}_{\min}\leq
\overline{K}_{\min}+(k-1)\overline{K}_{\max},$$ we have
\begin{eqnarray}\overline{K}_{\min}&\geq &\overline{Ric}^{(k)}_{\min}-(k-1)\overline{K}_{\max}\nonumber\\
&\geq & k\overline{Ric}^{(k)}_{\min}-(k-1)\overline{Ric}^{(k+1)}_{\max}.\end{eqnarray}
Suppose $\{e_1,e_2,e_3,e_4\}$ is an
orthonormal four-frame and $\lambda\in \mathbb{R}$.
From (2.1), (2.3), (4.8), (4.13) and (4.14) we get
\begin{eqnarray}&&R_{1313}+ R_{2323}-|R_{1234}|\nonumber\\
&\geq&
 \overline{Ric}^{(k)}_{\min}-\sum_{A= 3}^{k+1}\overline{R}_{A3A3}-\frac{2}{3}(\overline{K}_{\max}-\overline{K}_{\min})+\sum_{\alpha}\Big[h_{11}^{\alpha}h_{33}^{\alpha}+h_{22}^{\alpha}h_{33}^{\alpha}\nonumber\\
&&-\frac{3}{2}(h_{13}^{\alpha})^2
-\frac{3}{2}(h_{23}^{\alpha})^2-\frac{1}{2}(h^{\alpha}_{24})^2-\frac{1}{2}(h^{\alpha}_{14})^2\Big]\nonumber\\
&\geq& \overline{Ric}^{(k)}_{\min}-\Big(k-\frac{4}{3}\Big)[\overline{Ric}^{(k+1)}_{\max}-\overline{Ric}^{(k)}_{\min}]+\frac{2}{3}[k\overline{Ric}^{(k)}_{\min}-(k-1)\overline{Ric}^{(k+1)}_{\max}]\nonumber\\
&&+\sum_{\alpha}\Big[\sum_{i\neq
j}(h_{ij}^{\alpha})^2
 +\frac{(\sum_{i=1}^nh_{ii}^{\alpha})^2}{n-1}-S_{\alpha}
-\frac{3}{2}(h_{13}^{\alpha})^2
-\frac{3}{2}(h_{23}^{\alpha})^2-\frac{1}{2}(h^{\alpha}_{24})^2-\frac{1}{2}(h^{\alpha}_{14})^2\Big]\nonumber\\
&\geq&\frac{5k-1}{3}\Big[ \overline{Ric}^{(k)}_{\min}- \frac{5k-6}{5k-1} \overline{Ric}^{(k+1)}_{\max}\Big]+ \frac{n^2H^2}{n-1}-S.\end{eqnarray}
 Similarly, we have
\begin{equation}R_{1414}+ R_{2424}-|R_{1234}|\geq\frac{5k-1}{3}\Big[ \overline{Ric}^{(k)}_{\min}-\frac{5k-6}{5k-1} \overline{Ric}^{(k+1)}_{\max}\Big]+
\frac{n^2H^2}{n-1}-S.\end{equation}
 This together with (4.15) and the assumption implies
\begin{eqnarray}&&R_{1313}+\lambda^2 R_{1414} + R_{2323} +
\lambda^2R_{2424}- 2\lambda R_{1234}\nonumber\\
&\geq&R_{1313}+ R_{2323}-|R_{1234}|+\lambda^2(R_{1414}+ R_{2424}-|R_{1234}|)\nonumber\\
&\geq&(1+\lambda^2)\Big\{\frac{5k-1}{3}\Big[\overline{ Ric}^{(k)}_{\min}- \frac{5k-6}{5k-1} \overline{Ric}^{(k+1)}_{\max}\Big]+
\frac{n^2H^2}{n-1}-S\Big\}\nonumber\\
&>&0.
\end{eqnarray}
The assertion follows from Theorem A. \\
\hspace*{5mm}
(ii) Since \begin{equation}
\frac{\overline{Ric}^{(k)}_{\max}}{k}\geq \frac{\overline{Ric}^{(k+1)}_{\max}}{k+1},\end{equation}
we have$$\overline{Ric}^{(k)}_{\min}-\frac{(5k-6)(k+1)}{(5k-1)k}\overline{Ric}^{(k)}_{\max}\leq\overline{Ric}^{(k)}_{\min}-\frac{5k-6}{5k-1}\overline{Ric}^{(k+1)}_{\max}.$$
The assertion follows from the assumption and (i).\\
\hspace*{5mm}This completes the proof.\\\\
\hspace*{5mm}Taking $k=N-2$ in (i) of Lemma 4.3, we get the following theorem.\\\\
\textbf{Theorem 4.2.} \emph{Let $M^n$ be an $n(\geq4)$-dimensional
compact submanifold in an $N$-dimensional Riemannian manifold
$\overline{M}^{N}$. If
 $S<\frac{5N-11}{3}\Big[\overline{Ric}^{(N-2)}_{\min}-\frac{5N-16}{5N-11}\overline{Ric} _{\max}\Big]+
\frac{n^2H^{2}}{n-1}$, then $M$
is diffeomorphic to a spherical space form. In particular, if M is
simply connected, then M is diffeomorphic
to $S^n$.}\\\\
\textbf{Lemma 4.4.} \emph{Let $M^n$ be an $n(\geq4)$-dimensional
compact submanifold in an $N$-dimensional Riemannian manifold
$\overline{M}^{N}$. Suppose that $M$ satisfies one of the following conditions:\\
$(i)$
$S<\frac{5k^2+9k-8}{6}\Big[\overline{Ric}^{(k)}_{\min}-\frac{5k-6}{5k^2+9k-8}\overline{R}^{(k+1)}_{\max}\Big]+
\frac{n^2H^{2}}{n-1},$
 for some integer $k\in[2,N-2]$;\\$(ii)$ $S<\frac{5k^2+9k-8}{6}\Big[\overline{Ric}^{(k)}_{\min}-\frac{(k+2)(5k-6)}{s(5k^2+9k-8)}\overline{R}^{(k+1,s)}_{\max}\Big]+
\frac{n^2H^{2}}{n-1},$
 for some integers $k\in[2,N-2]$ and $s\in[2,k+1].$\\
  Then $M$ is diffeomorphic to a spherical space form. }\\\\
\textbf{Proof.} (i) It's seen from (2.8) and (2.9) that \begin{equation}\overline{K}_{\max}\leq \frac{1}{2}[\overline{R}^{(k+1)}_{\max}-(k+3)\overline{Ric}^{(k)}_{\min}],\end{equation} and $$\overline{Ric}^{(k)}_{\min}\leq
\overline{K}_{\min}+(k-1)\overline{K}_{\max},$$ which implies
\begin{eqnarray}K_{\min}&\geq &\overline{Ric}^{(k)}_{\min}-(k-1)\overline{K}_{\max}\nonumber\\
&\geq&
\frac{1}{2}[(k^2+2k-1)\overline{Ric}^{(k)}_{\min}-(k-1)\overline{R}^{(k+1)}_{\max}]
.\end{eqnarray} Suppose $\{e_1,e_2,e_3,e_4\}$ is an orthonormal
four-frame and $\lambda\in \mathbb{R}$. It follows from (2.1),
(2.3), (4.8), (4.19) and (4.20) that
\begin{eqnarray}&&R_{1313}+ R_{2323}-|R_{1234}|\nonumber\\
&\geq&
\overline{Ric}^{(k)}_{\min}-\sum_{A= 3}^{k+1}\overline{R}_{A3A3}-\frac{2}{3}(\overline{K}_{\max}-\overline{K}_{\min})+\sum_{\alpha}\Big[h_{11}^{\alpha}h_{33}^{\alpha}+h_{22}^{\alpha}h_{33}^{\alpha}\nonumber\\
&&-\frac{3}{2}(h_{13}^{\alpha})^2
-\frac{3}{2}(h_{23}^{\alpha})^2-\frac{1}{2}(h^{\alpha}_{24})^2-\frac{1}{2}(h^{\alpha}_{14})^2\Big]\nonumber\\
&\geq&\overline{Ric}^{(k)}_{\min}-\Big(\frac{k}{2}-
\frac{2}{3}\Big)[\overline{R}^{(k+1)}_{\max}-(k+3)\overline{Ric}^{(k)}_{\min}]\nonumber\\
&&
+\frac{1}{3}[(k^2+2k-1)\overline{Ric}^{(k)}_{\min}-(k-1)\overline{R}^{(k+1)}_{\max}]+\frac{n^2H^2}{n-1}-S\nonumber\\
&\geq&\frac{5k^2+9k-8}{6}\Big[\overline{Ric}^{(k)}_{\min}-\frac{5k
-6}{5k^2+9k-8}\overline{R}^{(k+1)}_{\max}\Big]+\frac{n^2H^2}{n-1}-S.\end{eqnarray}
 By a similar computation, we get
\begin{eqnarray}&&R_{1414}+ R_{2424}-|R_{1234}|\nonumber\\
&\geq&\frac{5k^2+9k-8}{6}\Big[\overline{Ric}^{(k)}_{\min}-\frac{5k
-6}{5k^2+9k-8}\overline{R}^{(k+1)}_{\max}\Big]+\frac{n^2H^2}{n-1}-S.\end{eqnarray}
 From (4.21) and (4.22), we obtain
\begin{eqnarray}&&R_{1313}+\lambda^2 R_{1414} + R_{2323} +
\lambda^2R_{2424}- 2\lambda R_{1234}\nonumber\\
&\geq&R_{1313}+ R_{2323}-|R_{1234}|+\lambda^2(R_{1414}+ R_{2424}-|R_{1234}|)\nonumber\\
&\geq& (1+\lambda^2)\Big\{\frac{5k^2+9k-8}{6} \Big[\overline{Ric}^{(k)}_{\min}-\frac{5k
-6}{5k^2+9k-8}\overline{R}^{(k+1)}_{\max}\Big]+\frac{n^2H^2}{n-1}-S\Big\}.\end{eqnarray}
 From (4.33), Theorem A and the assumption, we see that $M$ is diffeomorphic a spherical space form. \\
\hspace*{5mm} (ii) It follows from (2.9) that \begin{equation}
\frac{\overline{R}^{(k+1,s)}_{\max}}{s(k+1)}\geq\frac{\overline{R}^{(k+1
)}_{\max}}{(k+1)(k+2)},\end{equation} which implies
\begin{eqnarray}\overline{Ric}^{(k)}_{\min}-\frac{5k-6}{5k^2+9k-8}\overline{R}^{(k+1 )}_{\max}\geq\overline{Ric}^{(k)}_{\min}-\frac{(k+2)(5k-6)}{s(5k^2+9k-8)}\overline{R}^{(k+1,s)}_{\max} . \end{eqnarray}
Thus, the assertion follows from (i) and the assumption.\\ \hspace*{5mm}This completes the proof.\\\\
\textbf{Theorem 4.3.} \emph{Let $M^n$ be an $n(\geq4)$-dimensional
compact submanifold in an $N$-dimensional Riemannian manifold
$\overline{M}^{N}$. Assume that $M$ satisfies one of the following conditions:\\
$(i)$
$S<\frac{5N}{6s}\Big[\overline{Ric}_{\min}^{[s]}-\frac{s(5N^2-5N-12)}{5N}\overline{K}_{\max}\Big]+
\frac{n^2H^{2}}{n-1}; $\\
$(ii)$ $S< \frac{N(N+1)}{3s}\Big[\overline{Ric}_{\min}^{[s]}-\frac{s(N^2+2N+3)}{N(N+1)}\overline{Ric}_{\max}^{(N-2)}\Big]
+\frac{n^2H^2}{n-1}; $\\
 $(iii)$ $S< \frac{sN-s+6}{3}\Big [\overline{K}_{\min}-\frac{1}{sN-s+6}\overline{Ric}_{\max}^{[s]}\Big]+
\frac{n^2H^{2}}{n-1};$\\
$(iv)$ $S< \frac{5N^2-11N-6}{6}\Big[\overline{Ric}^{(N-2)}_{\min}-\frac{N(5N-16)}{s(5N^2-11N-6)}\overline{Ric}_{\max}^{[s]}\Big]+
\frac{n^2H^{2}}{n-1}$,\\
for some integer $s\in[2, N]$. Then the
normalized Ricci flow with initial metric $g_0$
$$\frac{\partial}{\partial t}g(t) = -2Ric_{g(t)} +\frac{2}{n}
r_{g(t)}g(t),$$ exists for all time and converges to a constant
curvature metric as $t\rightarrow\infty$. Moreover, $M$ is diffeomorphic to a spherical space form. In particular, if M is simply connected, then M is diffeomorphic
to $S^n$.}\\\\
 \textbf{Proof.} (i) Taking $k=N-1$ in (ii) of Lemma 4.2, we get the conclusion.\\
\hspace*{5mm}(ii) Since \begin{eqnarray} \overline{K}_{\min}&\geq &\frac{1}{2}[\overline{R}-(N+1)\overline{Ric}^{(N-2)}_{\max}]\nonumber\\
&\geq&\frac{1}{2}\Big[\frac{N\overline{Ric}_{\min}^{[s]}}{s}-(N+1)\overline{Ric}^{(N-2)}_{\max}\Big],\end{eqnarray}
we obtain\begin{eqnarray}\overline{K}_{\max}&\leq& \overline{Ric}^{(N-2)}_{\max}-(N-3)\overline{K}_{\min}\nonumber\\
&\leq&\frac{N^2-2N-1}{2}\overline{Ric}_{\max}^{(N-2)}-\frac{N(N-3)}{2s}\overline{Ric}_{\min}^{[s]}.\end{eqnarray}
Suppose $\{e_1,e_2,e_3,e_4\}$ is an orthonormal four-frame and
$\lambda\in \mathbb{R}$. Combing (2.1), (2.3), (4.8), (4.26) and
(4.27), we obtain
\begin{eqnarray}&&R_{1313}+ R_{2323}-|R_{1234}|\nonumber\\
&\geq&
2\overline{K}_{\min} -\frac{2}{3}(\overline{K}_{\max}-\overline{K}_{\min})+\sum_{\alpha}\Big[h_{11}^{\alpha}h_{33}^{\alpha}+h_{22}^{\alpha}h_{33}^{\alpha}\nonumber\\
&&-\frac{3}{2}(h_{13}^{\alpha})^2
-\frac{3}{2}(h_{23}^{\alpha})^2-\frac{1}{2}(h^{\alpha}_{24})^2-\frac{1}{2}(h^{\alpha}_{14})^2\Big]\nonumber\\
&\geq&\frac{n^2H^2}{n-1}-S+\frac{4}{3}\Big[\frac{N}{s}\overline{Ric}_{\min}^{[s]}-(N+1)\overline{Ric}^{(N-2)}_{\max}\Big]\nonumber\\
&&-\frac{1}{3}\Big[(N^2-2N-1)\overline{Ric}_{\max}^{(N-2)}-\frac{N(N-3)}{s}\overline{Ric}_{\min}^{[s]}\Big]\nonumber\\
&\geq&\frac{N+1}{3}\Big[\frac{N}{s}\overline{Ric}_{\min}^{[s]}-\frac{N^2+2N+3}{N+1}\overline{Ric}_{\max}^{(N-2)}\Big]+\frac{n^2H^2}{n-1}-S.\end{eqnarray}
 A similar discussion implies that
\begin{eqnarray}&&R_{1414}+ R_{2424}-|R_{1234}|\nonumber\\
&\geq&\frac{N+1}{3}\Big[\frac{N}{s}\overline{Ric}_{\min}^{[s]}-\frac{N^2+2N+3}{N+1}\overline{Ric}_{\max}^{(N-2)}\Big]+\frac{n^2H^2}{n-1}-S.\end{eqnarray}
 From (4.28) and (4.29), we get
\begin{eqnarray}&&R_{1313}+\lambda^2 R_{1414} + R_{2323} +
\lambda^2R_{2424}- 2\lambda R_{1234}\nonumber\\
&\geq&R_{1313}+ R_{2323}-|R_{1234}|+\lambda^2(R_{1414}+ R_{2424}-|R_{1234}|)\nonumber\\
&\geq& (1+\lambda^2)\Big\{\frac{N+1}{3}\Big[\frac{N}{s}\overline{Ric}_{\min}^{[s]}-\frac{N^2+2N+3}{N+1}\overline{Ric}_{\max}^{(N-2)}\Big]
 +\frac{n^2H^2}{n-1}-S\Big\}.\end{eqnarray}
Hence we get the conclusion from (4.30), the assumpion and Theorem A.\\
  \hspace*{5mm}(iii) It's seen from (2.9) that \begin{equation}\overline{K}_{\max}\leq \frac{1}{2}[\overline{Ric}^{[s]}_{\max}-(sN -s-2)\overline{K}_{\min}].\end{equation}Suppose $\{e_1,e_2,e_3,e_4\}$ is an
orthonormal four-frame and $\lambda\in \mathbb{R}$. By (2.1), (2.3),
(4.8), and (4.31), we have
\begin{eqnarray}&&R_{1313}+ R_{2323}-|R_{1234}|\nonumber\\
&\geq&
2\overline{K}_{\min} -\frac{2}{3}(\overline{K}_{\max}-\overline{K}_{\min})+\sum_{\alpha}\Big[h_{11}^{\alpha}h_{33}^{\alpha}+h_{22}^{\alpha}h_{33}^{\alpha}\nonumber\\
&&-\frac{3}{2}(h_{13}^{\alpha})^2
-\frac{3}{2}(h_{23}^{\alpha})^2-\frac{1}{2}(h^{\alpha}_{24})^2-\frac{1}{2}(h^{\alpha}_{14})^2\Big]\nonumber\\
&\geq&\frac{8}{3}\overline{K}_{\min}-\frac{1}{3}[\overline{Ric}^{[s]}_{\max}-(sN -s-2)\overline{K}_{\min}] +\frac{n^2H^2}{n-1}-S\nonumber\\
&\geq&\frac{1}{3}[(sN-s+6)\overline{K}_{\min}-\overline{Ric}^{[s]}_{\max}]+\frac{n^2H^2}{n-1}-S.\end{eqnarray}
 By using a similar argument, we get
\begin{eqnarray}&&R_{1414}+ R_{2424}-|R_{1234}|\nonumber\\
&\geq&\frac{1}{3}[(sN-s+6)\overline{K}_{\min}-\overline{Ric}^{[s]}_{\max}]++\frac{n^2H^2}{n-1}-S.\end{eqnarray}
 It follows from (4.32) and (4.33) that
\begin{eqnarray}&&R_{1313}+\lambda^2 R_{1414} + R_{2323} +
\lambda^2R_{2424}- 2\lambda R_{1234}\nonumber\\
&\geq&R_{1313}+ R_{2323}-|R_{1234}|+\lambda^2(R_{1414}+ R_{2424}-|R_{1234}|)\nonumber\\
&\geq& (1+\lambda^2)\Big\{\frac{1}{3}[(sN-s+6)\overline{K}_{\min}-\overline{Ric}^{[s]}_{\max}]+\frac{n^2H^2}{n-1}-S\Big\},\end{eqnarray}
which together with the assumption and Theorem A implies the conclusion.\\
\hspace*{5mm}(iv) The assertion follows from (ii) of Lemma 4.4 by taking $k=N-2$.\\
\hspace*{5mm}This proves the theorem.\\\\
\textbf{Theorem 4.4.} \emph{Let $M$ be a 3-dimensional compact
submanifold in an $N$-dimensional Riemannian manifold
$\overline{M}^N$. Assume that $M$ satisfies one of the following conditions:\\
$(i)$ $S<\frac{1}{2}[\overline{R}-(N^2-N-4)\overline{K}_{\max}]+
\frac{9}{2}H^2;$\\
$(ii)$ $S<\overline{R}-(N+1)\overline{Ric}^{(N-2)}_{\max}+\frac{9}{2}H^2;$\\
$(iii)$ $S<2\overline{K}_{\min}+\frac{9}{2}H^2;$\\
$(iv)$ $S<\frac{N^2-3N-2}{2} \Big[\overline{Ric}_{\min}^{(N-2)}-\frac{N-4}{N^2-3N-2}\overline{R}\Big]+\frac{9}{2}H^2$ for $N\geq 4$.\\
  Then $M$ is diffeomorphic to a spherical space form. In particular, if M is simply connected, then M is diffeomorphic
to $S^n$.}\\\\
\textbf{Proof.} For any unit tangent vector $u\in U_xM$
at $x\in M$, we choose an orthonormal three-frame $\{e_1,e_2,e_3\}$
such that $e_3=u$. \\
\hspace*{5mm}(i) From (2.1), (4.8) and the assumption, we obtain
\begin{eqnarray}Ric(u)&=&R_{1313}+ R_{2323}\nonumber\\
&\geq&
 \frac{1}{2}\Big(\overline{R}-2\sum_{A<B\neq 3}\overline{R}_{ABAB}\Big) +\sum_{\alpha}[h_{11}^{\alpha}h_{33}^{\alpha}+h_{22}^{\alpha}h_{33}^{\alpha}-(h_{13}^{\alpha})^2
-(h_{23}^{\alpha})^2 ]\nonumber\\
&\geq& \frac{1}{2}[\overline{R}-(N^2-N-4)\overline{K}_{\max}]\nonumber\\
&& +\sum_{\alpha}\Big[\sum_{i\neq
j}(h_{ij}^{\alpha})^2+\frac{(\sum_{i=1}^3h_{ii}^{\alpha})^2}{2}-S_{\alpha}
- (h_{13}^{\alpha})^2
-(h_{23}^{\alpha})^2 \Big]\nonumber\\
&\geq& \frac{1}{2}[\overline{R}-(N^2-N-4)\overline{K}_{\max}]+
\frac{9}{2}H^{2}-S\nonumber\\
&>&0. \end{eqnarray}
\hspace*{5mm}(ii) It follows from (2.8) that \begin{equation} \overline{K}_{\min}\geq \frac{1}{2}[\overline{R}-(N+1)\overline{Ric}^{(N-2)}_{\max}],\end{equation}
 which together with (2.1), (4.8) and the assumption implies
\begin{eqnarray}Ric(u)&=&R_{1313}+ R_{2323}\nonumber\\
&\geq&
2\overline{K}_{\min}+\sum_{\alpha}[h_{11}^{\alpha}h_{33}^{\alpha}+h_{22}^{\alpha}h_{33}^{\alpha}-(h_{13}^{\alpha})^2
- (h_{23}^{\alpha})^2]\nonumber\\
&\geq& \overline{R}-(N+1)\overline{Ric}^{(N-2)}_{\max}+\frac{n^2H^2}{n-1}-S\nonumber\\
&>&0.\end{eqnarray}
\hspace*{5mm}(iii)
 From (2.1), (4.8) and the assumption, we obtain
\begin{eqnarray}Ric(u)&=&R_{1313}+ R_{2323}\nonumber\\
&\geq&
2\overline{K}_{\min}+\sum_{\alpha}[h_{11}^{\alpha}h_{33}^{\alpha}+h_{22}^{\alpha}h_{33}^{\alpha} - (h_{13}^{\alpha})^2
- (h_{23}^{\alpha})^2]\nonumber\\
&\geq& 2\overline{K}_{\min}+\frac{n^2H^2}{n-1}-S\nonumber\\
&>&0.\end{eqnarray}
\hspace*{5mm}(iv) It's seen from (2.8) that \begin{equation} \overline{K}_{\max}\leq \frac{1}{2}[\overline{R}-(N+1)\overline{Ric}^{(N-2)}_{\min}].\end{equation}
 This together with (2.1), (4.8), (4.39)  and the assumption implies that
\begin{eqnarray}Ric(u)&=&R_{1313}+ R_{2323}\nonumber\\
&\geq&
\overline{Ric}_{\min}^{(N-2)}-(N-4)\overline{K}_{\max}+\sum_{\alpha}[h_{11}^{\alpha}h_{33}^{\alpha}+h_{22}^{\alpha}h_{33}^{\alpha} - (h_{13}^{\alpha})^2
- (h_{23}^{\alpha})^2]\nonumber\\
&\geq& \overline{Ric}_{\min}^{(N-2)}-\frac{N-4}{2}[\overline{R}-(N+1)\overline{Ric}^{(N-2)}_{\min}]+\frac{n^2H^2}{n-1}-S\nonumber\\
&\geq&\frac{N^2-3N-2}{2} \Big[\overline{Ric}_{\min}^{(N-2)}-\frac{N-4}{N^2-3N-2}\overline{R}\Big]+\frac{n^2H^2}{n-1}-S\nonumber\\
&>&0.\end{eqnarray}\hspace*{5mm} The assertion follows from (4.35), (4.37), (4.38), (4.40) and
Hamilton's theorem \cite{Hamilton}. This completes the proof of Theorem 4.4.\\\\
\textbf{Theorem 4.5.} \emph{Let $M^n$ be an $n(\geq4)$-dimensional
compact submanifold in an $N$-dimensional Riemannian manifold
$\overline{M}^{N}$. Assume that $M$ satisfies one of the following conditions:\\
$(i)$
$S<\frac{5}{6}N(N-1)(\overline{R}_0-\sigma_N\overline{K}_{\max})+
\frac{n^2H^{2}}{n-1},$\\
$(ii)$ $S<
\frac{N(N^2-1)}{3(N-2)}\Big[(N-2)\overline{R}_0-\mu_N\overline{Ric}_{\max}^{(N-2)}\Big]
+\frac{n^2H^2}{n-1},$\\
 $(iii)$ $S< \frac{N^2-N+6}{3}(\overline{K}_{\min}-\eta_N\overline{R}_0)+
\frac{n^2H^{2}}{n-1},$\\
$(iv)$ $S<
\frac{5N^2-11N-6}{6}\Big[\overline{Ric}^{(N-2)}_{\min}-\tau_N(N-2)\overline{R}_0\Big]+
\frac{n^2H^{2}}{n-1},$\\
where $\sigma_N$, $\mu_N$, $\eta_N$ and $\tau_N$  are defined as in
Theorems 1.1 and 1.2. Then $M$ is diffeomorphic to a spherical space form. In
particular, if M is simply connected, then M is diffeomorphic
to $S^n$.}\\\
\textbf{Proof.} By taking $k=N-1$ in (i) of Lemma 4.2, $s=N$ in (ii)
and (iii) of Theorem 4.3, and $k=N-2$ in (i) of Lemma 4.4,
respectively, we conclude that $M$ is diffeomorphic a spherical
space form. In particular, if $M$ is simply connected, then $M$ is
diffeomorphic
to $S^n$. This proves the theorem.\\\\
\section{Submanifolds
with weakly pinched curvatures} \hspace*{5mm}In this section, we
improve the differentiable sphere theorems \cite{Xu2} for
submanifolds with strictly pinched curvatures and obtain a
classification theorem for submanifolds with
weakly pinched curvatures.\\\\
\textbf{Proof of Theorem 1.3.}
If $n\geq4$, suppose $\{e_1,e_2,e_3,e_4\}$ is an
orthonormal four-frame and $\lambda\in \mathbb{R}$. From (2.1), (2.3) and (4.8), we have
 \begin{eqnarray}
 &&R_{1313}+\lambda^2 R_{1414} +\mu^2 R_{2323} +
\lambda^2\mu^2R_{2424}- 2\lambda \mu R_{1234}\nonumber\\
&\geq&(1+\lambda^2+\mu^2+\lambda^2\mu^2)\overline{K}_{\min}-\frac{4|\lambda\mu|}{3}(\overline{K}_{\max}-\overline{K}_{\min}) \nonumber\\
&&+\sum_{\alpha}\{h^{\alpha}_{11}h^{\alpha}_{33}-(h^{\alpha}_{13})^2+\mu^2\lambda^2\Big[h^{\alpha}_{22}h^{\alpha}_{44}-(h^{\alpha}_{24})^2]+2\mu\lambda h^{\alpha}_{14}h^{\alpha}_{23}\nonumber\\\
&&+ \mu^2[h^{\alpha}_{22}h^{\alpha}_{33}-(h^{\alpha}_{23})^2]+\lambda^2[h^{\alpha}_{11}h^{\alpha}_{44}-(h^{\alpha}_{14})^2]-2\mu\lambda h^{\alpha}_{13}h^{\alpha}_{24}\}\nonumber\\
&\geq&(1+\lambda^2+\mu^2+\lambda^2\mu^2)\frac{(4\overline{K}_{\min}-\overline{K}_{\max})}{3}\nonumber\\
&&+\sum_{\alpha}\Big\{\sum_{i<
j}(h_{ij}^{\alpha})^2+\frac{(\sum_{i=1}^nh_{ii}^{\alpha})^2}{2(n-1)}-\frac{S_{\alpha}}{2}-(h^{\alpha}_{13})^2-(h^{\alpha}_{14})^2-(h^{\alpha}_{23})^2\nonumber\\
&&+\mu^2\lambda^2\Big[\sum_{i<
j}(h_{ij}^{\alpha})^2+\frac{(\sum_{i=1}^nh_{ii}^{\alpha})^2}{2(n-1)}-\frac{S_{\alpha}}{2}-(h^{\alpha}_{24})^2-(h^{\alpha}_{14})^2-(h^{\alpha}_{23})^2\Big]\nonumber\\
&&+\mu^2\Big[\sum_{i<
j}(h_{ij}^{\alpha})^2+\frac{(\sum_{i=1}^nh_{ii}^{\alpha})^2}{2(n-1)}-\frac{S_{\alpha}}{2}-(h^{\alpha}_{23})^2-(h^{\alpha}_{13})^2-(h^{\alpha}_{24})^2\Big]\nonumber\\
&&+\lambda^2\Big[\sum_{i<
j}(h_{ij}^{\alpha})^2+\frac{(\sum_{i=1}^nh_{ii}^{\alpha})^2}{2(n-1)}-\frac{S_{\alpha}}{2}-(h^{\alpha}_{14})^2-(h^{\alpha}_{13})^2-(h^{\alpha}_{24})^2\Big]\Big\}\nonumber\\
&\geq&(1+\mu^2+\lambda^2+\lambda^2\mu^2)\Big[\frac{(4\overline{K}_{\min}-\overline{K}_{\max})}{3}+\frac{n^2H^2}{2(n-1)}-\frac{S}{2}\Big].
\end{eqnarray}
From (5.1) and the assumption that
$S\leq\frac{8}{3}\Big(\overline{K}_{\min}-\frac{1}{4}\overline{K}_{\max}\Big)+
\frac{n^2H^{2}}{n-1}$, we have
\begin{equation}R_{1313}+\lambda^2 R_{1414} +\mu^2 R_{2323} +
\lambda^2\mu^2R_{2424}- 2\lambda \mu R_{1234}\geq0,\end{equation}
i.e., $M\times \mathbb{R}^2$ has nonnegative isotropic curvature.\\
\hspace*{5mm}On the other hand, it follows from (2.1), (2.3) and (4.8) that\begin{eqnarray}&&R_{1313}+ R_{2323}-|R_{1234}|\nonumber\\
&\geq&
2\overline{K}_{\min}-\frac{2}{3}(\overline{K}_{\max}-\overline{K}_{\min})+\sum_{\alpha}\Big[h_{11}^{\alpha}h_{33}^{\alpha}+h_{22}^{\alpha}h_{33}^{\alpha}\nonumber\\
&&-\frac{3}{2}(h_{13}^{\alpha})^2
-\frac{3}{2}(a_{23}^{\alpha})^2-\frac{1}{2}(h^{\alpha}_{24})^2-\frac{1}{2}(h^{\alpha}_{14})^2\Big]\nonumber\\
 &\geq&\frac{8}{3}\Big(\overline{K}_{\min}-\frac{1}{4}\overline{K}_{\max}\Big)+ \frac{n^2H^2}{n-1}-S. \end{eqnarray}
The equalities in (5.3) hold only if
$$h_{ij}^{\alpha}=h_{33}^{\alpha}=0, \,\, \mbox{for all dinstinct}\,\, i,j \,\, \mbox{and any}\,\,  \alpha,$$and \begin{equation}h_{ii}^{\alpha}=h_{jj}^{\alpha},
 \,\, \mbox{for}\,\,  i,j\neq3 \,\, \mbox{and any}\,\, \alpha,\end{equation}
 i.e.,
\begin{equation}S=\frac{n^2H^2}{n-1}.\end{equation}
It follows from (5.3), (5.5) and the assumption that the Ricci curvature of $M$ is quasi-positive. This together with Aubin's theorem \cite{Aubin} implies that $M$ admits a metric with positive Ricci curvature.\\
\hspace*{5mm} By a similar discussion, we have
\begin{equation}R_{1414}+ R_{2424}-|R_{1234}|\geq \frac{8}{3}\Big(\overline{K}_{\min}-\frac{1}{4}\overline{K}_{\max}\Big)+
\frac{n^2H^2}{n-1}-S.\end{equation} The equality in (5.6) holds only if
$$h_{ij}^{\alpha}=h_{44}^{\alpha}=0, \,\, \mbox{for all dinstinct}\,\, i,j \,\, \mbox{and any}\,\,  \alpha,$$and \begin{equation}h_{ii}^{\alpha}=h_{jj}^{\alpha},
 \,\, \mbox{for}\,\,  i,j\neq4 \,\, \mbox{and any}\,\, \alpha.\end{equation} From (5.3), (5.5) and
(5.6), we have\begin{eqnarray}&&R_{1313}+R_{1414} + R_{2323} +
 R_{2424}- 2R_{1234}\nonumber\\
&=&(R_{1313}+ R_{2323}-|R_{1234}|)+ (R_{1414}+ R_{2424}-|R_{1234}|)\nonumber\\
&\geq&2\Big[\frac{8}{3}\Big(\overline{K}_{\min}-\frac{1}{4}\overline{K}_{\max}\Big)+
\frac{n^2H^2}{n-1}-S\Big],
\end{eqnarray}
  and the equality holds only if $S=\frac{n^2H^2}{n-1}$.  This together with the assumption and
  Lemma 2.3 implies that $M$ admits a metric with positive isotropic curvature. By Lemma 2.4, we conclude that $M$ is diffeomorphic to a spherical space form. In particular, if $M$ is simply connected, then $M$ is diffeomorphic to $S^n$.\\
  \hspace*{5mm}If $n=3$, for any unit tangent vector $u\in U_xM$
at $x\in M$, we choose an orthonormal three-frame $\{e_1,e_2,e_3\}$
such that $e_3=u$. From (4.8) and (4.38), we obtain$$Ric(u)=R_{1313}+ R_{2323} \geq 2\overline{K}_{\min}+\frac{n^2H^2}{n-1}-S,$$
 and the equality holds only if $S=\frac{n^2H^2}{n-1}$. Then from the assumption we know that $M$ has quasi-positive  Ricci curvature. Hence $M$ admits a metric with positive Ricci curvature by Aubin' theorem \cite{Aubin}.  This together
with Hamilton's theorem \cite{Hamilton} implies that $M$ is diffeomorphic to a spherical space form. In particular, if $M$ is simply connected, then $M$ is diffeomorphic to $S^n$.\\
\hspace*{5mm}This completes the proof of Theorem 1.3.\\\\
  \textbf{Corollary 5.1.} \emph{Let $M^n$ be an $n(\geq3)$-dimensional
complete submanifold in an $N$-dimensional Riemannian manifold
$\overline{M}^{N}$. If
$S\leq\frac{8}{3}\Big(\overline{K}_{\min}-\frac{1}{4}\overline{K}_{\max}\Big)+
\frac{n^2H^{2}}{n-1}$ and the strict inequality holds for some point
$x_{0}\in M$, then $M$ is diffeomorphic to
a spherical space form or $\mathbb{R}^n$. In particular, if M is simply connected, then M is diffeomorphic to $S^n$ or $\mathbb{R}^n$.}\\\\
\textbf{Proof.} From the assumption and Lemma 4.1 in \cite{Xu2}, we
know that $M$ has quasi-positive sectional curvature.
 When $M$ is noncompact, it follows from the Cheeger-Gromoll-Meyer-Perelman soul
 theorem \cite{Cheeger,Gromoll3,Perelman}
 that $M$ is diffeomorphic $\mathbb{R}^n$.
 When $M$ is compact, the assertion follows from Theorem 1.3. This  proves the corollary.\\\\
 \hspace*{5mm}For submanifolds in a sphere, we have the following theorem.\\\\
  \textbf{Theorem 5.1.} \emph{Let $M^n$ be an $n$-dimensional
 compact submanifold in the unit sphere $S^{N}$. Assume that $$S\leq2+\frac{n^2H^{2}}{n-1}.$$ We have the following possibilities:\\
 $(i)$ If $n=2$, then either $M$ is diffeomorphic  to $S^2$, $\mathbb{R}P^2$, or $M$ is flat.\\
 $(ii)$ If $n=3$, then M is diffeomorphic to a spherical space form.\\
 $(iii)$ If $n\geq 4$, then M is diffeomorphic to $S^n$.} \\\\
  \textbf{Proof.} If $n=2$, it's seen from the Gauss equation that $2K_M=2+4H^2-S$. This together with the assumption and Gauss-Bonnet theorem implies the conclusion.\\
  \hspace*{5mm}If $n=3$, we see from Proposition 2.1 in \cite{Shiohama2} that
\begin{eqnarray*}Ric(u)&\ge&\frac{2}{3}\Big[3+6H^2-S-\frac{3}{\sqrt{6}}H(S-3H^2)^{1/2}\Big]\\
&=&\frac{2}{3}\Big[(3+\frac{27}{4}H^2-\frac{3}{2}S)+\frac{3}{4}H^2+\frac{1}{2}(S-3H^2)
-\frac{3}{\sqrt{6}}H(S-3H^2)^{1/2}\Big]\\
&\geq&2+\frac{9}{2}H^2-S\end{eqnarray*} holds for any unit vector
$u\in T_x{M}$ at each point $x\in M,$ and the last inequality
becomes equality only if $S=\frac{9}{2}H^2.$
This together with the assumption implies $M$ has positive Ricci curvature. Hence the assertion follows from Hamilton's theorem \cite{Hamilton}.\\
\hspace*{5mm}If $n\geq 4$, from (4.8) we
get\begin{eqnarray} &&\sum_{k=2}^{n}[2|h(e_{1},e_{k})|^{2}-\langle
h(e_{1},e_{1}),h(e_{k},e_{k})\rangle]-(n-1)\nonumber \\
&=&\sum_{\alpha}\sum_{k=2}^n[2(h_{1k}^{\alpha})^2 -h_{11}^{\alpha}h_{kk}^{\alpha}]-(n-1)\nonumber\\
&\leq&\sum_{\alpha}\sum_{k=2}^n\Big\{2(h_{1k}^{\alpha})^2-\frac{1}{2}\Big[\sum_{i\neq
j}(h_{ij}^{\alpha})^2+\frac{(\sum_{i=1}^nh_{ii}^{\alpha})^2}{n-1}-S_{\alpha}\Big]\Big\}-(n-1)\nonumber\\
&\leq&\frac{n-1}{2}\Big(S-\frac{n^2H^2}{n-1}-2\Big).
\end{eqnarray}
The equalities in (5.9) hold only if $S=\frac{n^2H^2}{n-1}.$ From
the assumption $S\leq 2+\frac{n^2H^2}{n-1}$, we obtain
\begin{equation}\sum_{k=2}^{n}[2|h(e_{1},e_{k})|^{2}-\langle
h(e_{1},e_{1}),h(e_{k},e_{k})\rangle]-(n-1)<0.\end{equation} This
together with Theorem 2.1 implies that $M$ is simply connected. By
Theorem 1.3, we see that $M$ is diffeomorphic to $S^n$.
This proves Theorem 5.1.\\\\
  \hspace*{5mm}Moreover, we get the following  classification for complete submanifolds in an Euclidean space.\\\\
 \textbf{Theorem 5.2.} \emph{Let $M^n$ be an $n$-dimensional
oriented complete submanifold in the Euclidean space
$\mathbb{R}^{N}$.  Assume that
 $$S\leq\frac{n^2H^{2}}{n-1},\ \ H\neq 0.$$ We have the following possibilities:\\
  $(i)$ If $n=2$, then either $M$ is diffeomorphic to $S^2$, $\mathbb{R}^2$, or M is flat.\\
 $(ii)$ If $n=3$, then $M$ is either diffeomorphic to a spherical space form, $\mathbb{R}^3,$ or isometric to  $S^{2}(r_0)\times \mathbb{R}$.\\
   $(iii)$ If $n\geq 4$, then $M$ is either diffeomorphic $S^n$, $\mathbb{R}^n,$ or locally isometric to $S^{n-1}(r)\times \mathbb{R}$.}\\\\
\textbf{Proof.} It follows from the assumption and Lemma 2.1 that $K_M\geq 0$.\\
\hspace*{5mm} (i) Suppose that $M$ is compact. If $n=2$, it is seen from the Gauss-Bonnet theorem that $M$ is diffeomorphic to $S^2$ or $M$ is flat.\\
\hspace*{5mm}If $n=3$, we know that $Ric_M\geq 0$. This together
with Hamilton's theorem \cite{Hamilton} and
Lemma 2.1 implies that $M$ is diffeomorphic to a spherical space form, or $H$ is a
constant and $M$ is isometric to  $S^{2}\Big(\frac{2}{3H}\Big)\times \mathbb{R}$. Since $M$ is compact, the latter case is ruled out.\\
\hspace*{5mm}If $n\geq 4$, from the assumption
$S\leq\frac{n^2H^{2}}{n-1}$ and Lemma 2.1, we know $Ric_M\geq0$. We
claim that $M$ admits a metric with positive Ricci curvature.
Otherwise, it's seen from Aubin's theorem \cite{Aubin} that for each
point $x$ in $M$, there exists a unit vector $u$ in $T_{x}M$ such
that $Ric(u)=0$. By Lemma 2.1, we know that $H$ is constant and $M$
is isometric to $S^{n-1}\Big(\frac{n-1}{nH}\Big)\times\mathbb{ R}$,
which is noncompact. This contradicts to the compactness of $M$. By
the Bonnet-Myers theorem, the fundamental group $\pi_{1}(M)$ is
finite. Moreover, from the assumption, we know that
$S<\frac{n^2H^2}{n-2}$. It's seen from Lemma 2.2 and the universal
coefficient theorem that
 $H^{n-1}(M;\mathbb{Z})$ has no torsion, and hence neither does $H_{1}(M;\mathbb{Z})$
 by the Poincar$\acute{{\rm e}}$ duality. This together with the
 fact that
 $\pi_{1}(M)$ is finite implies $H_{1}(M;\mathbb{Z})=0$. Therefore we
 have $H_{n-1}(M;\mathbb{Z})=0$. Denote by $\widetilde{M}$ the universal Riemannian covering of
$M$. We may consider $\widetilde{M}$ be a Riemannian submanifold of
$\mathbb{R}^{N}$, and hence $\widetilde{M}$ is a homology sphere. Since
$\widetilde{M}$ is simply connected, it is a topological sphere, which
together with a result of Sjerve \cite{Sjerve} implies that $M$ is
simply connected. \\
\hspace*{5mm}On the other hand, from (5.1) and the assumption, we
know that $M\times \mathbb{R}^2$ has nonnegative isotropic
curvature. Moreover, it follows from (5.4), (5.7) and  the
assumption $H\neq 0$ that the equalities in (5.3) and (5.6) can not
hold simultaneously. Hence we see from (5.8) and the assumption that
$M$ has positive isotropic curvature. It follows from Lemma 2.4 that
  $M$ is diffeomorphic to a spherical space form. Since $M$ is simply connected, $M$ is diffeomorphic to $S^n$.\\
 \hspace*{5mm}(ii)  Suppose $M$ is noncompact. If $n=2$, it follows from the Cheeger-Gromoll-Meyer-Perelman soul theorem that $M$ is diffeomorphic to $\mathbb{R}^2$ or $M$ is flat.\\
   \hspace*{5mm}If $n=3$, a theorem due to Schoen-Yau \cite{Schoen1} and Zhu \cite{Zhu} states that if the Ricci curvature of $M$ is quasi-positive, then $M$ is diffeomorphic to $\mathbb{R}^3$. This together with Lemma 2.1 implies $M$ is diffeomorphic to $\mathbb{R}^3$ or isometric to $S^{2}(r_0)\times \mathbb{R} $.\\
    \hspace*{5mm}If $n\geq 4$,  we consider the following two
    cases:\\
     \hspace*{5mm}Case I. $K_M\geq 0$ and $K(\pi)>0$ for any 2-plane $\pi\subset T_{x_0}M$ at some point $x_0\in M$. By the
Cheeger-Gromoll-Meyer-Perelman soul theorem, we know that $M$ is
diffeomorphic to $\mathbb{R}^n$.\\
\hspace*{5mm}Case II. For each
point $x\in M$ there exists some 2-plane $\pi\subset T_xM$ such that
$K(\pi)=0.$ In this case, we get
$S\equiv\frac{n^2H^2}{n-1}$. Moreover, it follows from Lemma 2.1 and a result due to Ozgur(See \cite{Ozgur}, Theorem 4.1) that $M$ is conformally flat. This together with a theorem due to Carron \cite{Carron} implies that $M$ is diffeomorphic to $\mathbb{R}^n$ or locally isometric to $S^{n-1}(r)\times \mathbb{R}$.\\
 \hspace*{5mm} This completes the proof of Theorem 5.2.\\\\
\textbf{Proof of Theorem 1.4.} Combining Theorems 5.1 and 5.2 for $n\geq 4$, we complete the proof of Theorems 1.4.\\\\
\textbf{Corollary 5.2.} \emph{Let $M^n$ be an $n(\geq4)$-dimensional
oriented complete submanifold in an $N$-dimensional simply connected
space form $F^{N}(c)$ with $c\ge0$. Denote by $Ric^{[s]}(\cdot)$ the s-th weak Ricci curvature of M. Assume that
$$Ric^{[s]}_{\min}\ge
\frac{s(n+1)(n-2)c}{n}+\frac{sn(n-2)}{n-1}H^{2}$$ for some integer $s\in[1,n-1]$, where $c+H^2>0$. We have\\
$(i)$ If $c=0$, then M is either diffeomorphic $S^n$, $\mathbb{R}^n$, or locally isometric to $S^{n-1}(r)\times \mathbb{R}$.\\
$(ii)$ If M is compact, then M is diffeomorphic to $S^n$ .}\\\\
\textbf{Proof.} For $R\geq\frac{n}{s}Ric^{[s]}_{\min}$, the assertion follows from Theorem 1.4.\\\\
\textbf{\large{Acknowledgement.}}\ The authors would like to thank
Professors Kefeng Liu, Mu-Tao Wang, Shing-Tung Yau and Fangyang
Zheng for their
helpful discussions and valuable suggestions.\\\\

Juan-Ru Gu

Center of Mathematical Sciences\

Zhejiang University\

Hangzhou 310027\

China

E-mail address: gujr@cms.zju.edu.cn\\\\

Hong-Wei Xu

Center of Mathematical Sciences\

Zhejiang University\

Hangzhou 310027\

China

E-mail address: xuhw@cms.zju.edu.cn

\end{document}